\documentclass{article}

\usepackage{arxiv}
\usepackage{amssymb}
\usepackage{amsmath}
\usepackage{empheq}
\usepackage{amscd}
\usepackage{graphicx}
\usepackage{graphics}
\usepackage[noadjust]{cite}
\usepackage{amsthm}
\usepackage{caption}
\usepackage{multirow}
\usepackage{array}
\usepackage{rotating}

\usepackage{indentfirst}
\usepackage{tikz}
\usepackage{calc}
\usepackage{epsfig}
\usepackage{natbib}
\usepackage[makeroom]{cancel}
\usepackage{multicol}
\usepackage{csquotes}
\usepackage{algorithm}
\usepackage{algorithmic}
\usepackage{enumitem}
\usepackage{hyperref}
\usepackage{url}
\hypersetup{colorlinks,linkcolor={blue},citecolor={blue},urlcolor={blue}}

\usepackage{hyperref}       

\DeclareMathOperator*{\minimize}{minimize}
\newcommand{\vertiii}[1]{{\left\vert\kern-0.25ex\left\vert\kern-0.25ex\left\vert #1 
    \right\vert\kern-0.25ex\right\vert\kern-0.25ex\right\vert}}
    
\begin{document}
\title{Efficient extended-search space full-waveform inversion with unknown source signatures}

\author{\href{http://orcid.org/0000-0003-1805-1132}{\includegraphics[scale=0.06]{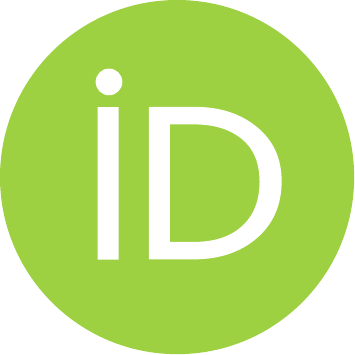}\hspace{1mm}Hossein S. Aghamiry} \\
  University Cote d'Azur - CNRS - IRD - OCA, Geoazur, Valbonne, France. 
  \texttt{aghamiry@geoazur.unice.fr}
  \And
\href{http://orcid.org/0000-0002-4981-4967}{\hspace{1mm}Frichnel W. Mamfoumbi-Ozoumet} \\ 
  University Cote d'Azur - CNRS - IRD - OCA, Geoazur, Valbonne, France. 
  \texttt{mamfoumbi@geoazur.unice.fr}
\And
  \href{https://orcid.org/0000-0002-9879-2944}{\includegraphics[scale=0.06]{orcid.pdf}\hspace{1mm}Ali Gholami} \\
  Institute of Geophysics, University of Tehran, Tehran, Iran.
  \texttt{agholami@ut.ac.ir} \\ 
  \And
\href{http://orcid.org/0000-0002-4981-4967}{\includegraphics[scale=0.06]{orcid.pdf}\hspace{1mm}St\'ephane Operto} \\ 
  University Cote d'Azur - CNRS - IRD - OCA, Geoazur, Valbonne, France. 
  \texttt{operto@geoazur.unice.fr}
  }


\renewcommand{\shorttitle}{IR-WRI with unknown sources, Aghamiry et al.}

\maketitle

\begin{abstract}
Full waveform inversion (FWI) requires an accurate estimation of source signatures. Due to the coupling between the source signatures and the subsurface model, small errors in the former can translate into large errors in the latter. 
When direct methods are used to solve the forward problem, classical frequency-domain FWI efficiently processes multiple sources for source signature and wavefield estimations once a single Lower-Upper (LU) decomposition of the wave-equation operator has been performed. However, this efficient FWI formulation is based on the exact solution of the wave equation and hence is highly sensitive to the inaccuracy of the velocity model due to the cycle skipping pathology. Recent extended-space FWI variants tackle this sensitivity issue through a relaxation of the wave equation combined with data assimilation, allowing the wavefields to closely match the data from the first inversion iteration. Then, the subsurface parameters are updated by minimizing the wave-equation violations. When the wavefields and the source signatures are jointly estimated with this approach, the extended wave equation operator becomes source dependent, hence making direct methods and, to a lesser extent, block iterative methods ineffective.
In this paper, we propose a simple method to bypass this issue and estimate source signatures efficiently during extended FWI. The proposed method replaces each source with a blended source during each data-assimilated wavefield reconstruction to make the extended wave equation operator source independent. Besides computational efficiency, the additional degrees of freedom introduced by spatially distributing the sources allows for a better signature estimation at the physical location when the velocity model is rough. 
We implement the source signature estimation with a variable projection method in the recently proposed iteratively-refined wavefield reconstruction inversion (IR-WRI) method.  Numerical tests on the Marmousi II and 2004 BP salt synthetic models confirm the efficiency and the robustness against velocity model errors of the new method compared to the case where source signatures are known. 
\end{abstract}
\section{Introduction}
%
%
%
%
Seismic wavefields carry information about subsurface and source, the latter being represented by its location and signature. In controlled source seismic on which this study is focused, the source locations are generally known accurately, while the source signatures are usually unknown and need to be estimated to perform reliable full-waveform inversion  (FWI) \citep{Tarantola_1984_ISR,Pratt_1998_GNF,Virieux_2009_OFW}. Furthermore, it is well acknowledged that the estimation of the source signature is easier in the frequency domain than in the time domain since the time-harmonic wave equation can be solved for each frequency separately \citep{Song_1995_FDAb,Pratt_1999_SWIb}. In frequency-domain seismic modeling, direct methods are the most suitable ones to process a large number of sources efficiently by forward/backward elimination, once a Lower-Upper (LU) decomposition of the so-called impedance matrix has been performed once \citep{Marfurt_1984_AFF}. When the size of the problem prevents using a direct solver, iterative methods speed-up the processing of multiple right-hand sides with block and recycling methods. Recently, \citet{Fang_2018_SEF} tackled the source signature estimation problem in the framework of an extended formulation of frequency-domain FWI called wavefield reconstruction inversion, where the recorded data are assimilated during the wavefield reconstruction to closely match the data with inaccurate subsurface models and hence prevent cycle skipping \citep{VanLeeuwen_2013_MLM}. In the formulation of \citet{Fang_2018_SEF}, the data assimilation makes the extended wave-equation operator source dependent, and hence the method is expensive since the LU decomposition needs to be performed for each source when a direct solver is used. Also, block-processing of multiple right-hand sides is not possible anymore with iterative methods. The focus of this paper is to revisit the source signature estimation problem in extended frequency-domain FWI such that the forward operator remains source-independent, and hence multiple sources can be processed efficiently as in the reduced-space FWI formulation.\\
%
%
%
%
Source signatures may be estimated before FWI or updated jointly with subsurface parameters during FWI iterations. For a fixed velocity model, the source signature estimation can be formulated as a least-squares quadratic data fitting problem \citep{Pratt_1999_SWIb}. The closed-form expression of the estimated source signature is given by the zero-lag cross-correlation between the calculated and recorded data, scaled by the auto-correlation of the calculated data. In this framework, the source signature estimation can be implemented in the classical reduced-space FWI iterations with two different approaches: In the first, the source signatures and the subsurface parameters are updated in an alternating mode, while the second approach enforces the closed-form expression of the estimated source signature as a function of the subsurface parameters in the objective function \citep{Aravkin_2012_ENP,Aravkin_2012_SEF,Liu_2013_AVP} through a variable projection approach \citep{Golub_2013_VPM}. \citet{Plessix_2011_GJI} review these two formulations in the frame of the adjoint-state method and conclude that the variable projection method is more versatile to implement the source signature estimation problem with specific data weighting, while \citet{Rickett_2013_VPM} showed that the variable projection approach was more resilient to phase errors in the wavelet than the alternated optimization. This source signature estimation doesn't introduce significant computational overhead in classical FWI since the gradient of the objective function with respect to the subsurface parameters is computed in the same way whether the source signature is available or estimated on the fly during the FWI iterations \citep{Aravkin_2012_ENP,Rickett_2013_VPM}.\\
%
%
%
%
%
%
In its more general form, FWI can be cast as a constrained optimization problem that aims to estimate the wavefields and the subsurface parameters by fitting the recorded data subject that the wave equation is satisfied \citep{Haber_2000_OTS}. Regardless of the source signature estimation issue, it is well acknowledged that FWI is highly nonlinear when the full search space encompassed by the wavefields and the subsurface parameters is projected onto the parameter space after elimination of the wavefield variables. This variable elimination, which is performed by forcing the wavefields to satisfy exactly the wave equation at each FWI iteration, makes FWI prone to cycle skipping as soon as the initial model is not accurate enough to predict recorded traveltimes with an error smaller than half a period \citep{Virieux_2009_OFW}. To mitigate the cycle skipping issue, some approaches implement the wave equation as a soft constraint with a penalty method such that the data can be closely matched with inaccurate subsurface models from the early FWI iterations \citep{Abubakar_2009_FDC,VanLeeuwen_2013_MLM,vanLeeuwen_2016_PMP}. Then, the subsurface model is updated by solving an overdetermined quadratic optimization problem, which consists of minimizing the source residuals generated by the wave equation relaxation. In these extended approaches, the wavefields are reconstructed by solving in a least-squares sense an overdetermined linear system gathering the wave equation weighted by the penalty parameter and the observation equation relating the simulated wavefield to the data through a sampling operator. In other words, the wavefields are reconstructed with data assimilation. This approach was called Wavefield Reconstruction Inversion (WRI) by \citet{VanLeeuwen_2013_MLM}. A variant of WRI, based upon the method of multipliers or augmented Lagrangian method,  was proposed by \citet{Aghamiry_2019_IWR} to increase the convergence rate and decrease the sensitivity of the algorithm to the relaxation (penalty) parameter choice. The augmented Lagrangian method combines a penalty method and a Lagrangian method, where the penalty term is used to implement the initial relaxation of the constraint while the Lagrangian term automatically tunes the sensitivity of the optimization to the constraint in iterations through the gradient ascent update of the Lagrange multipliers with the constraint violations. This method was called Iteratively-Refined(IR)-WRI, where the prefix IR refers to the iterative defect correction action of the Lagrange multipliers. \\
%
%
%
%
WRI was recently extended to jointly estimate the source signatures and the subsurface parameters \citep{Fang_2018_SEF}. In this approach, the monochromatic data-assimilated wavefield and the signature of the source are gathered in an unknown vector and are estimated in a least-squares sense. Although the wave-equation relaxation increases the robustness of the estimated source signatures against velocity model errors, the method is time-consuming because the normal system of the overdetermined wavefield-reconstruction problem is source-dependent, hence preventing efficient processing of multiple right-hand sides either with direct or iterative solvers. 
Another variant of WRI with unknown source signatures was proposed by \citet{Huang_2018_SEW} where WRI is re-parametrized in terms of extended sources and subsurface parameters. In this approach, the penalization (or annihilator) term is defined as the distance function from the real source position \citep[][ Their equation 9]{Huang_2018_SEW}, which means that the source signature estimation is implicitly embedded in the extended source reconstruction.
One issue with this approach is related to the presence of the Green functions in the Hessian of the extended source reconstruction subproblem, which makes the normal system very challenging to solve with a good accuracy \citep[][ Their equation 11]{Huang_2018_SEW}. Moreover, they update the subsurface parameters with a variable projection method, which precisely requires an accurate solution of the normal system for the extended sources.\\
%
%
%
%
In this paper, we implement a fast and robust (against model errors and recorded data inaccuracies) multi-source signature estimation in IR-WRI with a variable projection method, namely, the closed-form expression of the source signature is projected in the wavefield reconstruction subproblem.  To achieve the computational efficiency of the multi-source signature and wavefield reconstructions, we reconstruct each individual wavefield with blended sources. This blending makes the normal operator of the wavefield reconstruction subproblems source independent and hence amenable to efficient multi-source processing. Although we use source blending, we stress that we estimate one wavefield per physical source thanks to the assimilation of the source-dependent recorded data in the right-hand side of the normal system satisfied by the reconstructed wavefield. However, the source blending implies that, for each reconstructed wavefield, the source signature is a vector of dimension equal to the number of individual sources in the blended source.
When the velocity model is accurate, each entry of the signature vector is zero except the one located at the position of the physical source. When the velocity model is inaccurate, the other entries also contribute to decrease the data misfit, although their contribution is much less than the component located at the physical source position. Surprisingly, the additional degree of freedom provided by this spatially-distributed source signature helps to estimate a more accurate source signature when the velocity model is inaccurate, compared to the case where source blending is not used.  
The proposed algorithm solves IR-WRI with unknown source signatures in an alternating mode, first, jointly for data-assimilated wavefields and source signatures via variable projection, and then it solves a quadratic optimization for updating the model parameters. Because the extension created by the blended source assumption is artificial, we should correct its effects during the iteration. We propose two different algorithms for this correction.\\  
Numerical tests from the Marmousi II and 2004 BP salt model show that the proposed method is more robust (against inaccuracies in velocity model and data) and faster than traditional methods for source signature estimation and velocity model inversion. \\
%
This paper is organized as follows. In the method section, we first show how the source signature estimation can be combined with the extended wavefield reconstruction subproblem of IR-WRI by variable projection when each source is processed separately. We show that the extended wave equation operator becomes source dependent. The second part of the method section reviews the source blending approach that is used to make the extended wave equation operator source independent and hence amenable to efficient multi-source processing with direct methods. Two slightly different algorithms are proposed to implement the method. The paper continues with a numerical example section. We first assess the sensitivity of the source signature estimation to several parameters, such as the accuracy of the velocity model, noise, and the distance between sources and receivers. Then, we present applications of IR-WRI with the proposed efficient source signature estimation on the Marmousi II model and the BP salt model and compare the results when the source signatures are known and when they are estimated without the efficient blending strategy.
%
%
\section{Method}
\subsection{Notation and problem statement}
Frequency-domain FWI for multi-source acquisition with unknown source signatures can be formulated as the following constrained optimization problem \citep{Aghamiry_2019_IWR} 
\begin{equation} 
\minimize_{\bold{u}_i,{s}_i, \bold{m}\in \mathcal{M}}~~~~\mathcal{R}(\bold{m})  ~~~~\text{subject to}
~~~~ \begin{cases} \bold{A}(\bold{m})\bold{u}_i=s_i\boldsymbol{\phi}_i, ~~~i=1,2,...,n_s,\\ \bold{Pu}_i=\bold{d}_i, ~~~i=1,2,...,n_s \end{cases}
\label{main}
\end{equation} 
where $\bold{m} \in \mathbb{R}^{N\times 1}$ is the model parameter vector (squared slowness), $N$ is the number of discretized points of the medium, $n_s$ is the number of sources, $\bold{A}(\bold{m})=\Delta+\omega^2 \text{Diag}(\bold{m}) \in \mathbb{C}^{N\times N}$ is the Helmholtz operator, $\omega$ is the angular frequency, $\Delta$ is the Laplacian operator, Diag($\bold{x}$) denotes a diagonal matrix with the entries of the vector $\bold{x}$ on its diagonal, $\bold{u}_i \in \mathbb{C}^{N\times 1}$ and $\bold{d}_i \in \mathbb{C}^{M\times 1}$ denote the wavefield and the recorded data for the $i$'th source, respectively, $\bold{P} \in \mathbb{R}^{M\times N}$ is the observation operator and $M$ is the number of receivers. Also, ${s}_i \in \mathbb{C}$ is the source signature for the $i$'th source at frequency $\omega$, and $\boldsymbol{\phi}_i \in \mathbb{R}^{N \times 1}$ is a sparse vector defining the $i$'th source location. Finally, $\mathcal{R}(\bold{m})$ is an appropriate regularization function on the model domain and $\mathcal{M}$ is a convex set defined according to our prior knowledge of $\bold{m}$. For example, if we know the lower and upper bounds on $\bold{m}$ then
\begin{equation}
\mathcal{M} = \{\bold{m} \vert \bold{m}_{min} \leq \bold{m} \leq \bold{m}_{max}\}.
\end{equation}
IR-WRI solves the constrained problem \eqref{main} with the augmented Lagrangian method (or the method of multipliers). The augmented Lagrangian method combines a penalty term to relax the constraints during the early iterations and a Lagrangian term to control how accurately the constraint is satisfied at the convergence point \citep{Nocedal_2006_NO}. In this method, the primal variable and the Lagrange multiplier or dual variable are updated in alternating mode using a primal descent/dual ascent approach. Moreover, to make the computational cost tractable, we update the primal variables $\bold{u}$ and $\bold{m}$ in an alternating mode in the framework of the alternating-direction method of multipliers (ADMM) \citep{Boyd_2011_DOS}. The reader is referred to \citet{Aghamiry_2019_IWR}, \citet{Aghamiry_2019_IBC} and \citet{Aghamiry_2019_CRO} for more details about the ADMM-based IR-WRI algorithm. In the last two references, $\mathcal{R}(\bold{m})$ implements a total-variation (TV) regularization and an hybrid TV+Tikhonov regularization, respectively.
Compared to the above references, we extend IR-WRI to update the source signatures jointly with the wavefields during the wavefield reconstruction subproblem through a variable projection.\\
Beginning with an initial model $\bold{m}^{0}$ and assume $\bold{\hat{d}}_i^{0}=\bold{\hat{b}}^{0}_i=\bold{0}$, ADMM solves iteratively the multivariate optimization problem, equation \eqref{main}, with alternating directions as \citep[see][ for more details]{Boyd_2011_DOS,Aghamiry_2019_IWR}
\begin{subequations}
\label{ADMM}
 \begin{empheq}[left={\empheqlbrace\,}]{align}
(\bold{u}^{k+1}_i,s^{k+1}_i)&= \underset{\bold{u}_i,s_i}{\arg\min} ~ \Psi(\bold{u}_i,s_i,\bold{m}^{k},\bold{\hat{b}}_i^k,\bold{\hat{d}}_i^k),~~~~ i=1,2,...,n_s \label{primal_sig+wavefield}\\
\bold{m}^{k+1}&= \underset{\bold{m}\in \mathcal{M}}{\arg\min} ~ \sum_{i=1}^{n_s}\Psi(\bold{u}_i^{k+1},\bold{m},s^{k+1}_i,\bold{\hat{b}}_i^k,\bold{\hat{d}}_i^k) \label{primal_sigma}\\
\bold{\hat{b}}_i^{k+1} &= \bold{\hat{b}}_i^k  +{s}^{k+1}_i\boldsymbol{\phi}_i- \bold{A}(\bold{m}^{k+1})\bold{u}^{k+1}_i,~~~~i=1,2,...,n_s\label{dual_b}\\ 
\bold{\hat{d}}_i^{k+1} &= \bold{\hat{d}}_i^k  +\bold{d}_i- \bold{P}\bold{u}_i^{k+1},~~~~ i=1,2,...,n_s  \label{dual_d}
\end{empheq}
\end{subequations}
where
\begin{align}
\Psi(\bold{u}_i,s_i,\bold{m},\bold{\hat{b}}_i^k,\bold{\hat{b}}_i^k) =
\mathcal{R}(\bold{m}) + 
\|\bold{Pu}_i-\bold{d}_i-\bold{\hat{d}}_i^k\|_2^2 +\lambda \|\bold{A}(\bold{m})\bold{u}_i-s_i\boldsymbol{\phi}_i-\bold{\hat{b}}_i^k\|_2^2, \label{eqpsi}
\end{align}
is the scaled form of the augmented Lagrangian \citep[][ section 3.1.1]{Boyd_2011_DOS}, $\bullet^k$ is the value of $\bullet$ at iteration $k$, the scalar $\lambda>0$ is the penalty parameter assigned to the wave equation constraint, and $\bold{\hat{b}}_i^{k}\in \mathbb{C}^{N\times 1}$, $\bold{\hat{d}}^{k}_i\in \mathbb{C}^{M\times 1}$ are the scaled Lagrange multipliers, which are updated through a dual ascent scheme by the running sum of the constraint violations (source and data residuals) as shown by equations~\eqref{dual_b}-\eqref{dual_d}. The penalty parameter $\lambda$ can be tuned in equation \eqref{eqpsi} such that the estimated wavefields approximately fit the observed data from the first iteration at the expense of the accuracy with which the wave equation is satisfied, while the iterative update of the Lagrange multipliers progressively corrects the errors introduced by these penalizations such that both of the observation equation and the wave equation are satisfied at the convergence point with acceptable accuracies.\\
Here, we focus on the optimization subproblem \eqref{primal_sig+wavefield}. The readers are referred to \citet{Aghamiry_2019_IBC,Aghamiry_2019_CRO,Aghamiry_2020_FWI} for the closed-form expression of the optimization subproblem \eqref{primal_sigma} with bound constraints and different regularizations and to \citet{Aghamiry_2020_RWI} for a more robust implementation of this subproblem against velocity model errors with phase retrieval. \\
Optimization problem \eqref{primal_sig+wavefield} is quadratic in $\bold{u}_i$ and $s_i$ and can be written as
\begin{equation} \label{main_obj_single}
\underset{\substack{\bold{u}_i,{s}_i}}{\arg\min} ~~\|\bold{Pu}_i-\bold{d}_i\|_{{2}}^2 + \lambda \|\bold{A}_k\bold{u}_i-{s}_i\boldsymbol{\phi}_i\|_{{2}}^2,   
\end{equation}
where $\bold{A}_k \equiv \bold{A}(\bold{m}^k)$. \citet{Fang_2018_SEF} solve Eq. \eqref{main_obj_single} jointly for $\bold{u}_i$ and ${s}_i$ by gathering them in a single vector and solve an $(N+1)\times(N+1)$ linear system (Their Eq. 24) instead of the $N \times N$ system of the original WRI \citep{VanLeeuwen_2013_MLM}. Here, we want to solve Eq. \eqref{main_obj_single} with the variable projection method \citep{Golub_2013_VPM}.
By taking derivative of Eq. \eqref{main_obj_single} with respect to $s_i$, we get that ${s}_i$ satisfies $\boldsymbol{\phi}_i^{\text{T}}(\bold{A}_k\bold{u}_i-{s}_i\boldsymbol{\phi}_i)=0$, where $\bullet^{\text{T}}$ denotes the conjugate transpose of $\bullet$, and since $\boldsymbol{\phi}_i^{\text{T}}\boldsymbol{\phi}_i=1$, we get 
\begin{equation}
{s}_i=\boldsymbol{\phi}_i^{\text{T}}\bold{A}_k\bold{u}_i. 
\label{eqsi}
\end{equation}
Substituting the above expression of ${s}_i$ into Eq. \eqref{main_obj_single} leads to a mono-variate optimization problem for the wavefield as
\begin{equation} \label{main_obj_U_single}
\underset{\substack{\bold{u}_i}}{\arg\min} ~~\|\bold{Pu}_i-\bold{d}_i\|_{\text{2}}^2 + \lambda \|\bold{Q}_i\bold{A}_k\bold{u}_i\|_{\text{2}}^2,   
\end{equation}
where $\bold{Q}_i = \bold{I} - \boldsymbol{\phi}_i\boldsymbol{\phi}_i^{\text{T}}$, and $\bold{I}$ is the identity matrix.
In equation~\eqref{main_obj_U_single}, the matrix $\boldsymbol{\phi}_i\boldsymbol{\phi}_i^{\text{T}}$ is a diagonal matrix with one nonzero coefficient equal to 1 at the location of the source $i$, and $\bold{Q}_i$ is another diagonal matrix complementary to $\boldsymbol{\phi}_i\boldsymbol{\phi}_i^{\text{T}}$: its diagonal entries equal to 1 except at the source position where the coefficient is zero. 
The second term in equation~\eqref{main_obj_U_single} penalizes the predicted source $\bold{A}_k\bold{u}_i$ at all spatial points except at the physical source location consistently with the elimination (or projection) of $s_i$, equation~\eqref{eqsi}, from the optimization variables. 
Minimization of Eq. \eqref{main_obj_U_single} with respect to $\bold{u}_i$ gives (note that $\bold{Q}_i^{\text{T}}\bold{Q}_i=\bold{Q}_i^2=\bold{Q}_i$)
\begin{equation} \label{U_single}
\bold{u}_i^{k+1}=\left(\bold{P}^{\text{T}}\bold{P} + \lambda \bold{A}_k^{\text{T}}\bold{Q}_i\bold{A}_k\right)^{-1} \bold{P}^{\text{T}}\bold{d}_i.
\end{equation}
The explicit relation between the estimated source signature and the data can be obtained as
\begin{equation}
 \label{S_single}
{s}_i^{k+1}=\boldsymbol{\phi}_i^{\text{T}}\bold{A}_k\bold{u}_i^{k+1}=\boldsymbol{\phi}_i^{\text{T}}\bold{A}_k\left(\bold{P}^{\text{T}}\bold{P} + \lambda \bold{A}_k^{\text{T}}\bold{Q}_i\bold{A}_k\right)^{-1} \bold{P}^{\text{T}}\bold{d}_i =\boldsymbol{\phi}_i^{\text{T}}\left(\bold{G}_k^{\text{T}}\bold{G}_k + \lambda \bold{Q}_i\right)^{-1} \bold{G}_k^{\text{T}}\bold{d}_i,
\end{equation}
where $\bold{G}_k=\bold{P}\bold{A}_k^{-1}$ is the rank-deficient forward operator sampling the Green function $\bold{A}_k^{-1}$ at receiver positions. 
Equation \eqref{S_single} shows that the source signature is estimated by first back propagating the data in time from the receiver positions, i.e. $\bold{G}_k^{\text{T}}\bold{d}_i$, and then correct the blurring effects induced by the limited bandwidth of the data and the limited spread of the receivers by applying the inverse of the Hessian, i.e. $\left(\bold{G}_k^{\text{T}}\bold{G}_k + \lambda \bold{Q}_i\right)$.\\
The optimization problem, Eq. \eqref{main_obj_U_single} and its closed-form solution, Eq. \eqref{U_single}, share some similarities with the extended source reconstruction method described in \citet[their Eqs. 10 and 11]{Huang_2018_SEW} as a source-independent variant of WRI, although there are two differences in their formulation: first, their state variables are the extended sources instead of the extended wavefields, i.e., $\bold{b}^e_i=\bold{A}_k\bold{u}_i$ where $\bold{b}^e_i$ is the $i$'th extended source; Second, they used another anihilator function for $\bold{Q}_i$ (Their equation 9), which is zero at the source location and linearly increases away from it. The difficulty with the method of \cite{Huang_2018_SEW} is related to the presence of the Green's functions $\bold{A}_k^{-1}$ in the Hessian $\left(\bold{G}_k^{\text{T}}\bold{G}_k + \lambda \bold{Q}_i\right)$ which makes the system very difficult, if not impossible, to solve exactly.\\
\subsection{Efficient multi-source processing with source blending}
The proposed method for joint estimation of source signature and data-assimilated wavefield, Eq. \eqref{main_obj_U_single}, is robust against velocity model error but it is computationally-expensive because the normal operator $\left(\bold{P}^{\text{T}}\bold{P} + \lambda \bold{A}_k^{\text{T}}\bold{Q}_i\bold{A}_k\right)$ is source-dependent, hence preventing efficient multi-source processing with either direct or iterative solvers. This difficulty also exists in the method of \citet{Fang_2018_SEF} for joint estimation of wavefield and source signature (Their Eq. 24). Although they proposed a block matrix formulation to overcome the computational challenges (Their Eq. 33), it seems inefficient and not applicable for large scale problems. \\
To overcome this computational issue and design a fast and accurate multi-source signature and wavefield estimation, we assume that a virtual blended source generates each wavefield.
The true signature of this blended source is the physical source signature at the physical source location and is zero elsewhere. 
By doing so, each source can be written as $\boldsymbol{\Phi}\bold{s}_i$ where $\bold{\Phi}=[\boldsymbol{\phi}_1~\boldsymbol{\phi}_2~...~\boldsymbol{\phi}_{n_s}] \in \mathbb{R}^{N \times n_s}$ is a tall matrix including the shifted delta functions ($\boldsymbol{\phi}_i$) in its columns and the true $\bold{s}_i \in \mathbb{C}^{n_s \times 1}$ is a vector that contains the $i$'th physical source signature ($s_i$) and it is zero elsewhere. We stress at this stage that this reformulation of the source, $\boldsymbol{\Phi}\bold{s}_i$, is equivalent to the original one $s_i \Phi_i$ given in equation~\eqref{main}. \\
Plugging the new source expression in the objective function \eqref{main_obj_single} and taking the derivative with respect to $\bold{s}_i$ gives
\begin{equation}
\bold{s}_i=\boldsymbol{\Phi}^T\bold{A}_k\bold{u}_i.
\label{eqsib}
\end{equation} 
Projecting this expression into the cost function and remembering that $\boldsymbol{\Phi}^T\boldsymbol{\Phi}=\bold{I}$ give 
\begin{equation} \label{main_obj_U_multi}
\underset{\substack{\bold{u}_i}}{\arg\min} ~~\|\bold{Pu}_i-\bold{d}_i\|_{\text{2}}^2 + \lambda \|\bold{Q}\bold{A}_k\bold{u}_i\|_{\text{2}}^2,   
\end{equation}
where $\bold{Q} = \bold{I} - \bold{\Phi}\bold{\Phi}^{\text{T}}=\frac{1}{n_s}\sum_{i=1}^{n_s} \bold{Q}_i$. \\
The closed-form expression of the wavefield $\bold{u}_i$ at iteration $k+1$ is now given by
\begin{equation} \label{U_multi}
\bold{u}_i^{k+1}=\left(\bold{P}^{\text{T}}\bold{P} + \lambda \bold{A}_k^{\text{T}}\bold{Q}\bold{A}_k\right)^{-1} \bold{P}^{\text{T}}\bold{d}_i.
\end{equation}
Comparing Eq. \eqref{U_multi} and Eq. \eqref{U_single} shows that the Hessian $\left(\bold{P}^{\text{T}}\bold{P} + \lambda \bold{A}_k^{\text{T}}\bold{Q}\bold{A}_k\right)$ is source independent, hence preserving the benefit of direct solver method to process efficiently multiple sources once one LU factorization has been performed. This solves the computational issue. However, the new parametrization of the source makes the optimization problem blind to the fact that the vector $\bold{s}_i \in \mathbb{C}^{n_s \times 1}$ should have only one non zero entry at index $i$. It gives equal probability to all source positions to reconstruct the wavefield $\bold{u}_i$, equation~ \eqref{main_obj_U_multi}, hence leading to a blended wavefield. This is highlighted by the fact that $\bold{A}_k\bold{u}_i$ is potentially dense in the closed-form expression of the reconstructed signature $\bold{s}_i$, equation~\eqref{eqsib}, when the velocity model is inaccurate. Conversely, when the velocity model is the true one, $\bold{A}_{true}\bold{u}_i= s_i \Phi_i$, is sparse.
%
%
%
%
%
%
The blended source assumption gives an extra degree of freedom to the optimization problem to decrease the cost function. Surprisingly, we will show in the \textit{Numerical results} section that the estimated source signature obtained with the source blending approach, Eq. \eqref{main_obj_U_multi}, is more accurate than the counterpart obtained without blending, Eq. \eqref{main_obj_U_single}, when we start IR-WRI from a rough initial velocity model.\\
The source blending is artificial and its effects (extra non zero coefficients in $\bold{s}_i$ vectors) must be removed during iterations. In the following, we propose two algorithms to achieve this goal. For the sake of compact notations, we recast from now the optimization problem, eq. \eqref{main}, in matrix form.\\
Multi sources can be processed efficiently in frequency-domain modeling by gathering them in the right-hand side (rhs) of the Helmholtz system in a matrix format. Considering $n_s$ sources, the multi-rhs Helmholtz system is written as $\bold{AU}=\boldsymbol{\Phi}\bold{S}$ where $\bold{U}=[\bold{u}_1~\bold{u}_2~ ...~ \bold{u}_{n_s}] \in \mathbb{C}^{N \times n_s}$ and $\bold{S} \in \mathbb{C}^{n_s \times n_s}$ is a square diagonal matrix with the source signatures on its main diagonal, $\bold{S}_{ii}=s_i$. Introducing the data matrix $\bold{D}=[\bold{d}_1~\bold{d}_2~ ...~ \bold{d}_{n_s}] \in \mathbb{C}^{M \times n_s}$ gives the following optimization problem for multi-source signature and wavefield estimations 
\begin{equation} \label{main_obj_F}
\underset{\substack{\bold{U},\bold{S}}}{\arg\min} ~~\|\bold{PU}-\bold{D}\|_{\text{F}}^2 + \lambda \|\bold{A}_k\bold{U}-\bold{\Phi}\bold{S}\|_{\text{F}}^2,   
\end{equation}
where $\|\bullet\|_{\text{F}}^2$ denotes the Frobenius norm. 
Solving Eq. \eqref{main_obj_F} with a variable projection method gives 
\begin{equation}
\bold{S}=\bold{\Phi}^{\text{T}}\bold{A}_k\bold{U},
\label{eqS}
\end{equation}
where the diagonal components of $\bold{S}$ are dominant, and the off-diagonal coefficients represent the non-physical source components associated with each of the $n_s$ blended sources. 
Plugging the explicit expression of $\bold{S}$ into Eq. \eqref{main_obj_F} leads to a mono-variate optimization problem for $\bold{U}$, the closed-form expression of which is given by
\begin{equation} \label{U}
\bold{U}^{k+1}=\left(\bold{P}^{\text{T}}\bold{P} + \lambda \bold{A}_k^{\text{T}}\bold{Q}\bold{A}_k\right)^{-1} \bold{P}^{\text{T}}\bold{D}.
\end{equation}
Eq. \eqref{U} is the same as Eq. \eqref{U_multi} but for multi-data $\bold{D}$. 
At this stage, we didn't impose any constraint on the structure of $\bold{S}$ which is potentially dense when the velocity model is inaccurate.
Plugging the expression of $\bold{U}^{k+1}$ in equation~\eqref{eqS} gives the explicit expression of $\bold{S}$ 
\begin{equation} \label{S}
\bold{S}^{k+1}=\boldsymbol{\Phi}^T\bold{A}_k\bold{U}^{k+1}=\boldsymbol{\Phi}^T\bold{A}_k\left(\bold{P}^{\text{T}}\bold{P} + \lambda \bold{A}_k^{\text{T}}\bold{Q}\bold{A}_k\right)^{-1} \bold{P}^{\text{T}}\bold{D}.
\end{equation}  
Before proceeding with the subsurface parameter updating, we extract the approximate signature of the physical sources as $\bold{s}=\text{diag}(\bold{S})$, where diag($\bold{X}$) denotes a vector that contains the diagonal elements of matrix $\bold{X}$. The different steps of IR-WRI with source signature estimation are reviewed in Algorithm \eqref{alg1}, which begins with an initial model $\bold{m}^{0}$ and initial dual variables $\bold{\hat{B}}^0=\bold{0}$ and $\bold{\hat{D}}^0=\bold{0}$.\\
\begin{algorithm}
 \caption{IR-WRI with unknown sources} \label{alg1}
 \begin{algorithmic}[1]
 \REQUIRE Initial velocity model $\bold{m}^0$
 \STATE set $\bold{\hat{B}}^0=\bold{0}$ and $\bold{\hat{D}}^0=\bold{0}$ 
 \REPEAT
 \STATE $\bold{U}^{k+1}= \underset{\substack{\bold{U}}}{\arg\min} ~ \|\bold{PU}-\bold{D}-\bold{{\hat{D}}}^k\|_{\text{F}}^2+\lambda \|\bold{QA}_{k}\bold{U} - \bold{\hat{B}}^{k}\|_{\text{F}}^2$ 
 \STATE $\bold{S}^{k+1} = \text{Diag}(\text{diag} (\bold{\Phi}^{\text{T}}\bold{A}_{k}\bold{U}^{k+1}))$ 
\STATE $\bold{m}^{k+1}= \underset{\bold{m}}{\arg\min} ~ \mathcal{R}(\bold{m})+\lambda \|\bold{A}(\bold{m})\bold{U}^{k+1}-\bold{S}^{k+1}-\bold{\hat{B}}^{k}\|_{\text{F}}^2 $
\STATE $\bold{{\hat{B}}}^{k+1} = \bold{\hat{B}}^k+\bold{S}^{k+1}- \bold{A}_{k+1}\bold{U}^{k+1}$ 
\STATE $\bold{\hat{D}}^{k+1} = \bold{\hat{D}}^k  +\bold{D}- \bold{P}\bold{U}^{k+1}$ 
\UNTIL{stopping conditions are satisfied.}
\end{algorithmic}
\end{algorithm}
The $\bold{U}$-subproblem in Algorithm \ref{alg1} introduce errors in the extended wavefield reconstruction due to the source blending. These errors can be corrected iteratively by the action of the Lagrange multipliers, which are formed by the source residuals computed without source blending (Line 6 of algorithm \ref{alg1}). This is shown by the fact that the diagonal components of the source signature matrix are used (Line 4 of algorithm \ref{alg1}) instead of the whole matrix, equation~\eqref{S}.  
The right-hand side correction term $\bold{\hat{B}}^k$ in the objective function of the $\bold{U}$-subproblem in Algorithm \ref{alg1} gathers the running sum of the source residuals of previous iterations. This iterative refinement leads to the error forgetting property discussed by \citet{Yin_2013_EFG} in the frame of Bregman iterations, which means that the error correction performed at the current iteration is made independent of the error corrections performed at previous iterations. Here, this iterative solution refinement by right-hand side updating is necessary to correct three sources of errors: the first results from the fact that each primal subproblem is solved keeping fixed the other primal variable, the second from the fact that we solve a constrained problem with a penalty method keeping the penalty parameter fixed and the third from the non-physical source blending.
Another application of iterative refinement in AVO inversion is presented in \citet{Gholami_2017_CNA}, where the linearized Zoeppritz equations are used to simplify the primal problem, while the dual problem compensates the linearization-related errors by computing the residuals with the exact Zoeppritz equations.\\
When we seek to reconstruct a complicated velocity model starting from a rough initial model, the Algorithm \ref{alg1} wasn't able to fully remove the detrimental effects of the source blending and failed to reach the same minimizer as IR-WRI with a known source. This prompts us to propose Algorithm \ref{alg2} that includes one extra step compared to Algorithm \ref{alg1} by re-estimating the wavefields (Line 5 of Algorithm \ref{alg2}) with the diagonal restriction of the source signature matrix (Line 4 of Algorithm \ref{alg2}). By doing so, the pollution effects of the non-physical sources are removed from the reconstructed wavefields. The improvement provided by this wavefield refinement is illustrated in the next \textit{Numerical results} section.
\begin{algorithm}[H]
 \begin{algorithmic}[1]
 \caption{IR-WRI with unknown sources with wavefield correction} \label{alg2}
 \REQUIRE starting point $\bold{m}^0$
 \STATE set $\bold{\hat{B}}^0=\bold{0}$ and $\bold{\hat{D}}^0=\bold{0}$ 
 \REPEAT
 \STATE $\bold{U}^{k+\frac12}= \underset{\substack{\bold{U}}}{\arg\min} ~ \|\bold{PU}-\bold{D}-\bold{{\hat{D}}}^k\|_{\text{F}}^2+\lambda \|\bold{QA}_{k}\bold{U} - \bold{\hat{B}}^{k}\|_{\text{F}}^2$ 
 \STATE $\bold{S}^{k+1} = \text{Diag}(\text{diag} (\bold{\Phi}^{\text{T}}\bold{A}_{k}\bold{U}^{k+\frac12}))$
 \STATE $\bold{U}^{k+1}= \underset{\substack{\bold{U}}}{\arg\min} ~ \|\bold{PU}-\bold{D}-\bold{{\hat{D}}}^k\|_{\text{F}}^2+\lambda \|\bold{A}_{k}\bold{U} - \bold{S}^{k+1}-\bold{\hat{B}}^{k}\|_{\text{F}}^2$   
\STATE $\bold{m}^{k+1}= \underset{\bold{m}}{\arg\min} ~ \mathcal{R}(\bold{m})+\lambda \|\bold{A}(\bold{m})\bold{U}^{k+1}-\bold{S}^{k+1}-\bold{\hat{B}}^{k}\|_{\text{F}}^2 $
\STATE $\bold{{\hat{B}}}^{k+1} = \bold{\hat{B}}^k+\bold{S}^{k+1}- \bold{A}_{k+1}\bold{U}^{k+1}$ 
\STATE $\bold{\hat{D}}^{k+1} = \bold{\hat{D}}^k  +\bold{D}- \bold{P}\bold{U}^{k+1}$ 
\UNTIL{stopping conditions are satisfied.}
\end{algorithmic}
\end{algorithm}
Algorithm \ref{alg2} requires two LU decompositions at each IR-WRI iteration (one for $\bold{U}^{k+\frac12}$ and the other one for $\bold{U}^{k+1}$), but still remains much cheaper than the algorithm performing one LU decomposition per source. \\
%
%
%
%
\section{Numerical results}
\begin{figure}[htb!]
\centering
\includegraphics[scale=1.3]{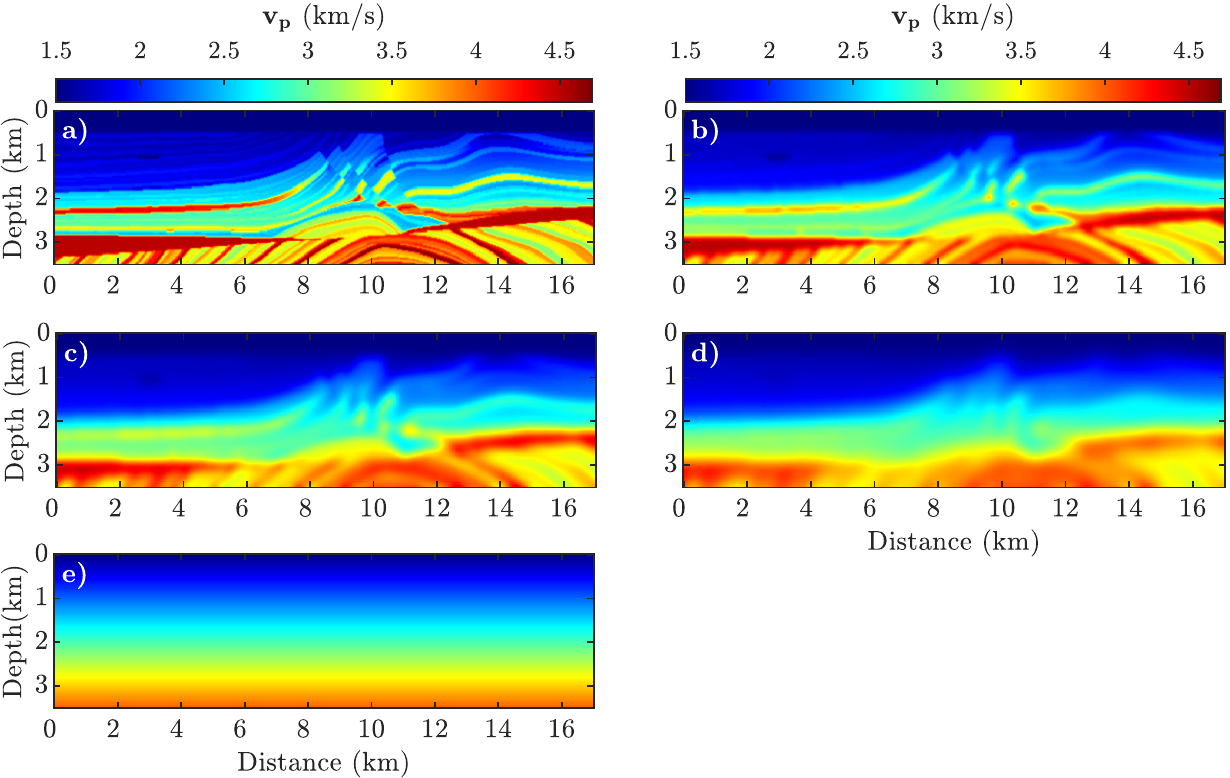}
\caption{Marmousi II model. (a) True model. (b-e) Smoothed versions of the true model, referred to as models 1 (b), 2 (c), 3 (d) and 4 (f) in the paper.}
\label{Fig_Mar}
\end{figure}

We first investigate different aspects of the proposed method for efficient source signature estimation in IR-WRI (referred to as {\it{joint}} approach in the following) with the Marmousi II model (Fig. \ref{Fig_Mar}a) and compare its performance when each source is processed separately in IR-WRI (referred to as {\it{separate}} approach in the following), equation~\eqref{ADMM}, and when the source signatures are estimated with the method of \citet{Pratt_1999_SWIb} (referred to as {\it{conventional}} approach),  as
 \begin{equation} \label{Pratt}
 \bold{S}=\left( \bold{\Phi}^T\bold{G}^T\bold{G}\bold{\Phi} \right)^{-1} \bold{\Phi}^T\bold{G}^T \bold{D}.
 \end{equation} 
Then, we compare the performances of the proposed Algorithms \ref{alg1} and \ref{alg2} with those of separate approaches when the source signatures are known (classical IR-WRI) and unknown, equation~\eqref{ADMM}. We use the Marmousi II model and a scaled version of the left target of the challenging 2004 BP salt model \citep{Billette_2004_BPB} for this comparison.\\   
For all the numerical tests, we use a 9-point finite-difference staggered-grid stencil with PML boundary condition (along the model's edges except for the top where the free-surface boundary condition is used) and anti-lumped mass to solve the Helmholtz equation.
\begin{figure}[htb!]
\centering
\includegraphics[scale=0.72]{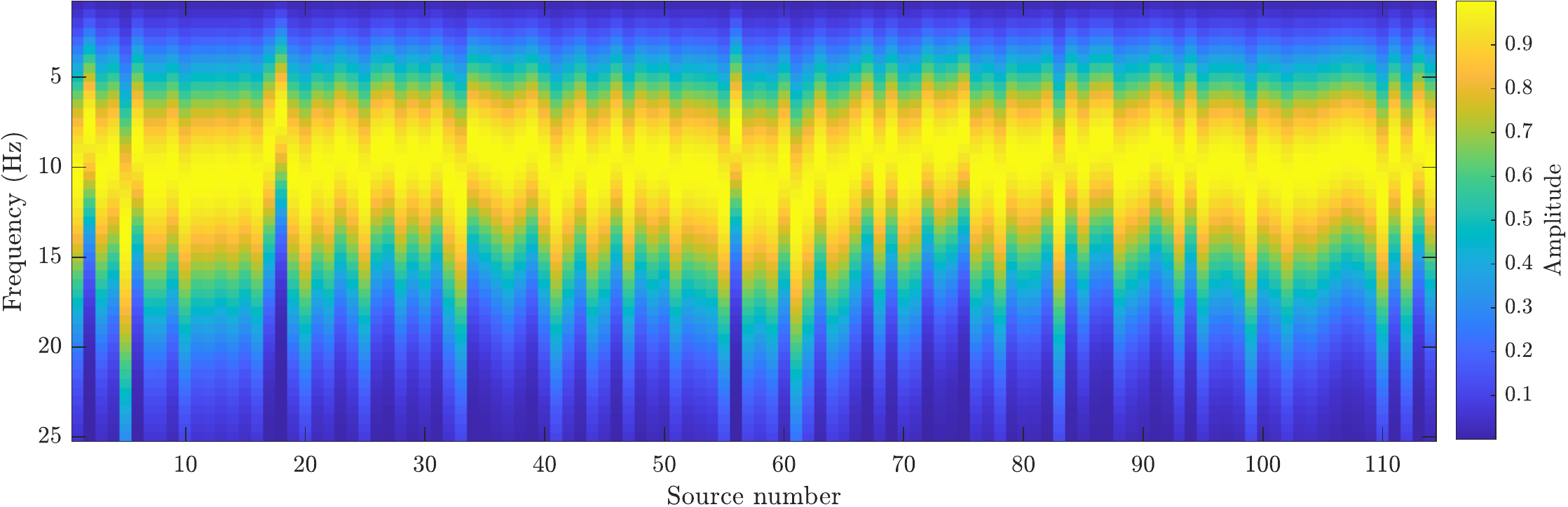}
\caption{Amplitude spectrum of source signatures used in Marmousi II test.}
\label{Fig_wavelet}
\end{figure}
\subsection{Marmousi II test}
We first illustrate the performance of the separate, joint and conventional methods for source signature estimation and investigate the robustness of these methods against the accuracy of the initial velocity model, noise in the recorded data, the vertical distance between the source and the receiver profiles, and the number of sources. We use four initial models for these tests (Fig. \ref{Fig_Mar}b-\ref{Fig_Mar}e), which are referred to as models 1
-4. 
%
%
The fixed-spread acquisition contains 114 point sources, the source signatures of which are Ricker wavelets of different central frequency and initial phase (Figure \ref{Fig_wavelet}), and a line of receivers spaced 50~m apart at the surface. The central frequency for each source signature is selected randomly between [7-15]~Hz, and the peak of each wavelet is centered randomly between [0-0.4]~s. 
\subsubsection{Sensitivity to the background velocity model}
\begin{figure}
\centering
\includegraphics[scale=0.75]{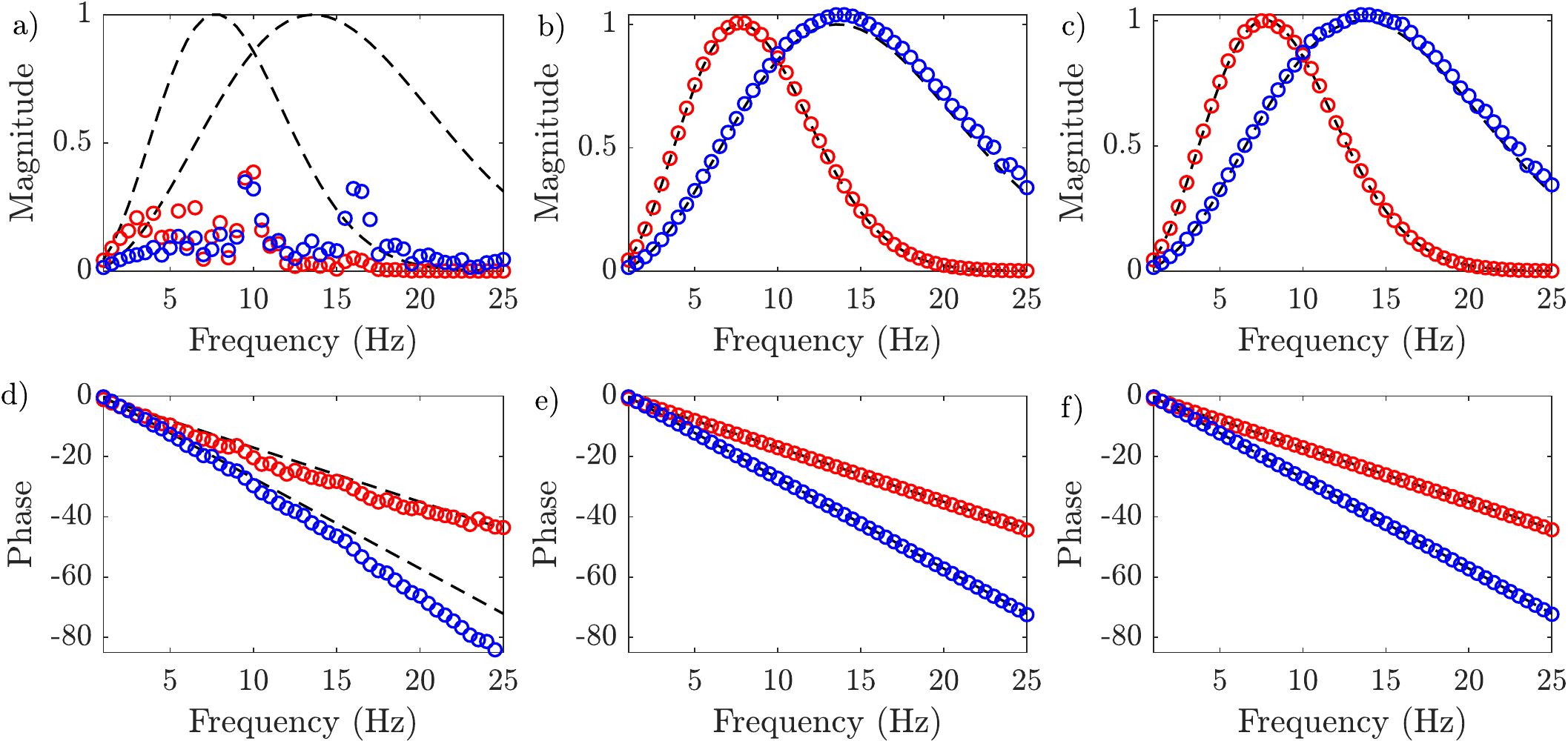}
\caption{Estimated source signature (magnitude and unwrapped phase) with different methods for a couple of sources with source numbers 18 and 33 in Fig. \ref{Fig_wavelet} when model 4 (Fig. \ref{Fig_Mar}e) is used as the initial velocity model. (a) The conventional method of \citet{Pratt_1999_SWIb}, (b) the separate method, and (c) the joint method. The blue and red points show the estimated source signatures for the first and second sources, respectively, and the black dashed lines show the true ones. (d-e) same as (a-c) but for unwrapped phase.}
\label{Fig_err_wavelet_test}
\end{figure}

\begin{figure}
\centering
\includegraphics[scale=0.5]{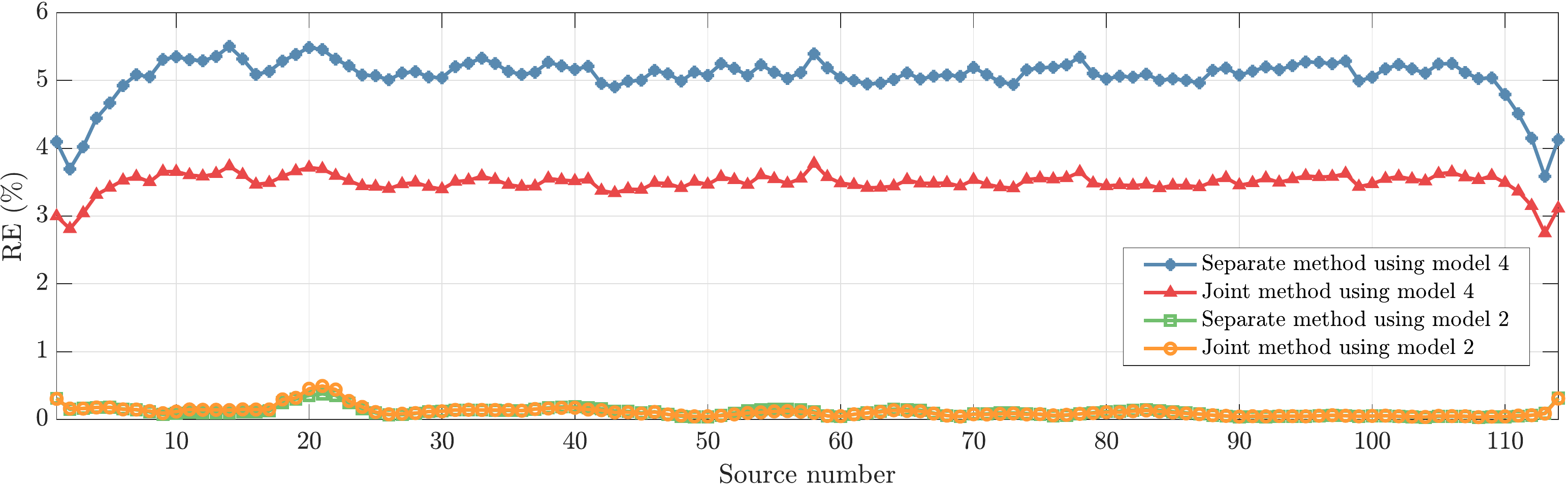}
\caption{RE of the estimated source signatures for separate and joint method when a rough initial model 4 (Fig. \ref{Fig_Mar}e) and the kinematically accurate initial model 2 (Fig. \ref{Fig_Mar}c) are used.}
\label{Fig_err_wavelet}
\end{figure}

First, we put the line of sources at 75~m depth, generate the data with the true velocity model (Fig. \ref{Fig_Mar}a) and use the 1D gradient velocity model (model 4, Fig. \ref{Fig_Mar}e) to estimate the source signatures with the conventional method  (Eq. \eqref{Pratt}), the separate method  (Eq. \eqref{S_single}), and the joint method (Eq. \eqref{S}).
We show the magnitude and phase of the estimated source signatures for a couple of sources with source numbers 18 and 33 (Fig. \ref{Fig_wavelet}) in Fig. \ref{Fig_err_wavelet_test}. 
First, the results clearly show the improvement achieved by the relaxed wave-equation methods (the separate and joint methods) (Fig.\ref{Fig_err_wavelet_test}b,c,e,f) compared to the conventional method (Fig.\ref{Fig_err_wavelet_test}a,d). Second, both separate and joint methods estimate accurate source signatures but with a different computational burden (one LU decomposition for the joint method against 114 LU decompositions for the separate method). \\
To gain more quantitative insights on the accuracy of methods, we plot the relative error (RE) of the estimated source signatures as a function of the source number for each method in Fig. \ref{Fig_err_wavelet}. 
The RE for the estimated source signature is defined as
\begin{equation}
\text{RE}= \frac{\sqrt{\sum_{j=f_1}^{f_n}[\bold{x}_j^*-\hat{\bold{x}}_j]^2}}{\sqrt{\sum_{j=f_1}^{f_n}{\bold{x}^*_{j}}^2}},
\end{equation}
where $\bold{x}^*_j$ and $\hat{\bold{x}}_j$ are the true and estimated source signatures at frequency $j$, respectively, and $f_1$ and $f_n$ are the minimum and maximum frequencies, respectively. 
In this figure, we show the RE of the estimated source signatures for the separate and joint methods when the rough velocity model 4 (Fig. \ref{Fig_Mar}e) and the kinematically accurate model 2 (Fig. \ref{Fig_Mar}c) are used as background velocity model. We don't show the RE of the conventional method in this figure because it is much higher than those obtained with the separate and joint approaches. It is shown that when the initial velocity model is rough, the joint method performs better than the separate method (the blue and red curves in Fig. \ref{Fig_err_wavelet}). This probably results from the extra degrees of freedom available in the joint method compared to the separate counterpart. More precisely, the error in the estimated source signature generated by the inaccuracy of the initial velocity model is entirely mapped at the physical location of the source in the separate method. In contrast, this error is distributed across the different components of the blended source in the joint method, hence, providing a better estimation of the source signature at the location of the physical source.
On the other hand, as the velocity model becomes more accurate, the separate and joint methods reach the same accuracy for the source signature estimation (the green and orange curves in Fig. \ref{Fig_err_wavelet}). \\
%
%
%
%
%
%
\subsubsection{Sensitivity to noise}
\begin{figure}
\centering
\includegraphics[scale=0.83]{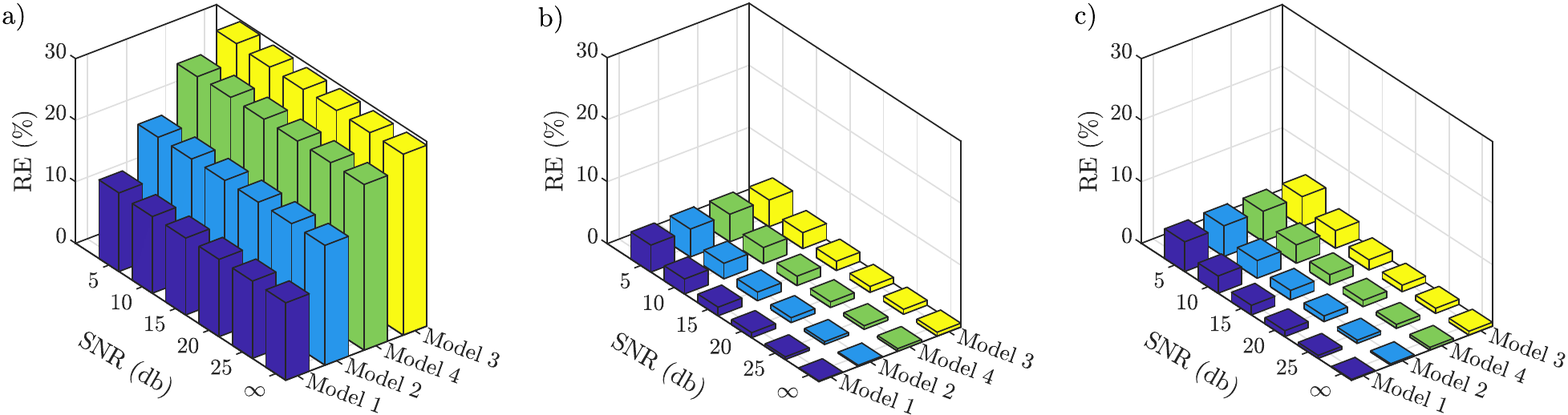}
\caption{Average RE over all the sources for the estimated source signatures of (a) conventional method, (b) separate method and (c) joint method as a function of SNR of recorded data and initial velocity model.}
\label{Fig_SNR}
\end{figure}

We continue by assessing the robustness of the methods against the noise in the recorded data and the error in the initial velocity model. We repeat the same test as before with different initial velocity models and different levels of random Gaussian noises in the data. The average RE over all the sources for the conventional, separate, and joint methods are shown in Figs. \ref{Fig_SNR}a-\ref{Fig_SNR}c, respectively, as functions of the signal to noise ratio (SNR) of the data and the initial velocity model. In this paper, the SNR of data is defined as 
\begin{equation}
\text{SNR}=20 \log \left( \frac{Clean~data~RMS~amplitude}{Noise~RMS~amplitude} \right).
\end{equation}
First of all, the conventional method is robust against the noise in the data, but it is sensitive to the errors in the initial velocity model, as illustrated in the previous section. In contrast, the separate and joint methods are robust against inaccuracy of the initial velocity model due to the extended search space allowing for data fitting with an inaccurate model but are sensitive to the noise in recorded data due to the risk of noise overfitting. 
Although the RE of the separate and joint methods increases with the amount of noise, it remains however far less than the conventional method. It can be seen that even in the worst scenario (with the lowest SNR data and a rough initial model), such methods based on wave-equation relaxation still can estimate an acceptable source signature.\\ 
\subsubsection{Sensitivity to the distance between the source and the receivers}
\begin{figure}
\centering
\includegraphics[scale=1]{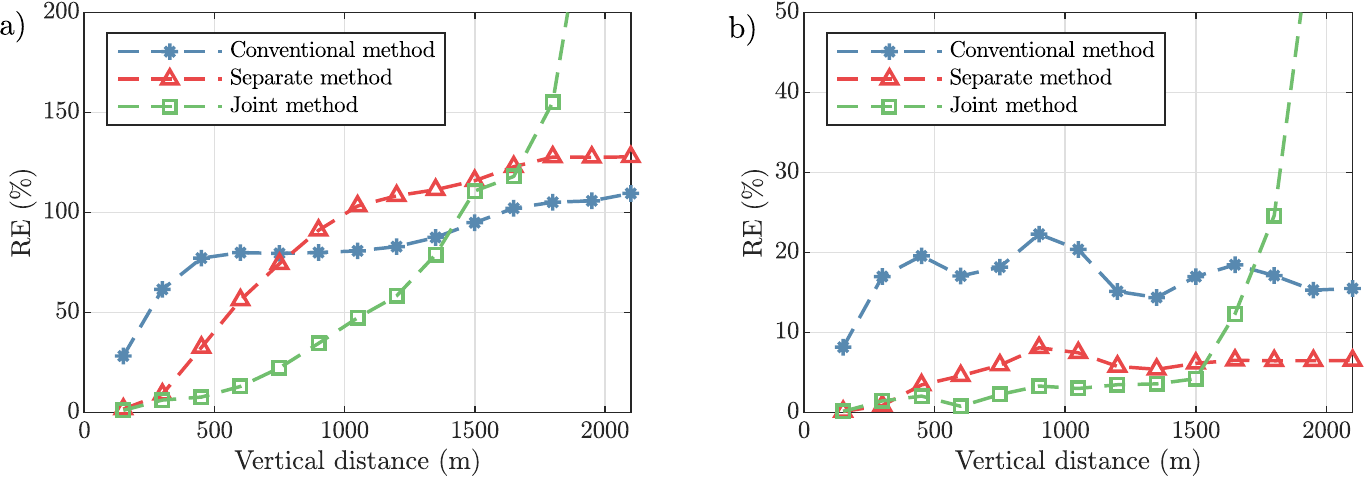}
\caption{Average RE over all the sources for the estimated source signatures as a function of the vertical distance between the source and receiver lines when (a) the rough velocity model 4 (Fig. \ref{Fig_Mar}e), (b) the kinematically accurate velocity model 2 (Fig. \ref{Fig_Mar}c) are used.}
\label{Fig_Depth}
\end{figure}

In IR-WRI, the accuracy of the estimated source signature is directly controlled by the distance between the source and the closest receiver. This results because the extended wavefield, from which the source signature is estimated, equation~\eqref{eqS}, matches well the recorded data only near the receivers when the background velocity model is inaccurate. This implies that when the source is close to a receiver, it will be estimated from an accurate estimation of the wavefield at this receiver. 
Moreover, the impedance matrix in equation~\eqref{eqS}, which is built from a potentially inaccurate velocity model, will not generate significant errors when applied to the wavefield to generate the source when the latter and the receiver are close to each other. In this section, we show the robustness of the different methods against this distance. We do the same test as before, but we change the depth of the source line while keeping the receiver line at the surface. The average RE of the estimated source signatures, summed over all the sources, are plotted in Fig. \ref{Fig_Depth} as a function of the depth of the source profile for the conventional, separate, and joint methods and for the rough model 4 (Fig. \ref{Fig_Depth}a) and the kinematically accurate model 2 (Fig. \ref{Fig_Depth}b). First, it can be seen that for the vertical distances up to 1500 m, the joint method has a better performance for both initial velocity models, but it becomes unstable beyond this vertical distance where the subsurface becomes more complex (green curves in Fig. \ref{Fig_Depth}). This suggests that the additional degrees of freedom in the joint method drives the least-squares problem, Eq. \eqref{main_obj_F}, toward an inaccurate local minimizer when the extended wavefield becomes too inaccurate at the source location. For surface acquisitions, this should however not be an issue in practice. For towed-streamer acquisitions, the sources are close to the nearest receiver, and both of them are in the water. In seabed acquisition, the reciprocal sources are on the seabed and potentially far away from the receivers. However, the medium between the receiver and the source layouts (i.e., the water) is known. On land, areal acquisitions are classically designed with sources and receivers at the surface with short nearest offset.
\subsubsection{Sensitivity to the number of sources}
\begin{figure}
\centering
\includegraphics[scale=0.8]{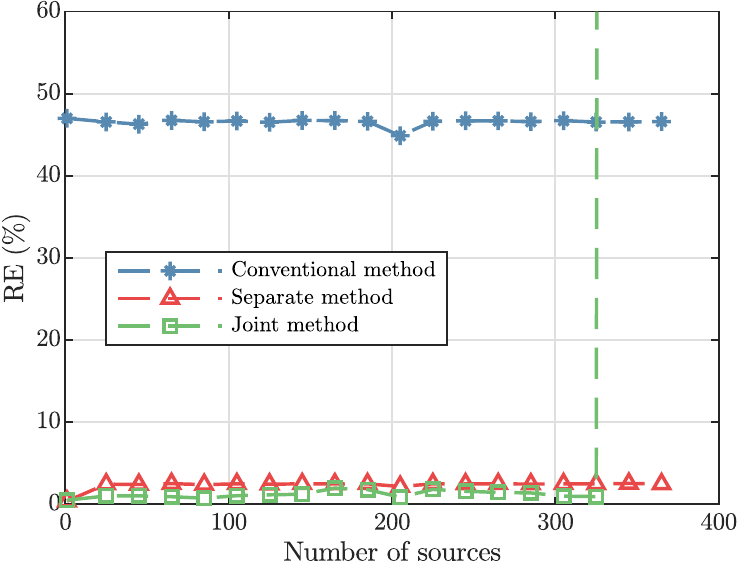}
\caption{Average RE over all the sources for the estimated source signatures as a function of the number of the sources when the rough initial velocity model 4 (Fig. \ref{Fig_Mar}e) and a line of 320 receivers at the surface are used.}
\label{Fig_NOS}
\end{figure}

The other aspect that we need to investigate is the number of sources and receivers. We repeat the same test as before several times with the rough initial velocity model 4 and a line of 320 receivers with 50~m spacing at the surface. 
For each of them, we use a line of sources with a different number of sources ranging from 1 to 360 with an interval 20 (sources) deployed at 75~m depth. The average RE over all the sources for the estimated source signatures as a function of the number of sources for the conventional, separate, and joint methods is shown in Fig. \ref{Fig_NOS}. Both the conventional method and the separate method have a stable behavior (blue and red curves), while the joint method becomes unstable as soon as the number of sources exceeds the number of receivers. 
In this case, the source signature estimation problem becomes under-determined in the case of the joint method, and equation Eq. \eqref{main_obj_F} converges to the least-norm solution. 
In practice, this issue can be easily bypassed by subdividing the sources into patches of suitable dimension, gathering possible closely-spaced sources. \\ 
\subsubsection{Assessment of IR-WRI with source signature estimation}
\begin{figure}
\centering
\includegraphics[scale=0.95]{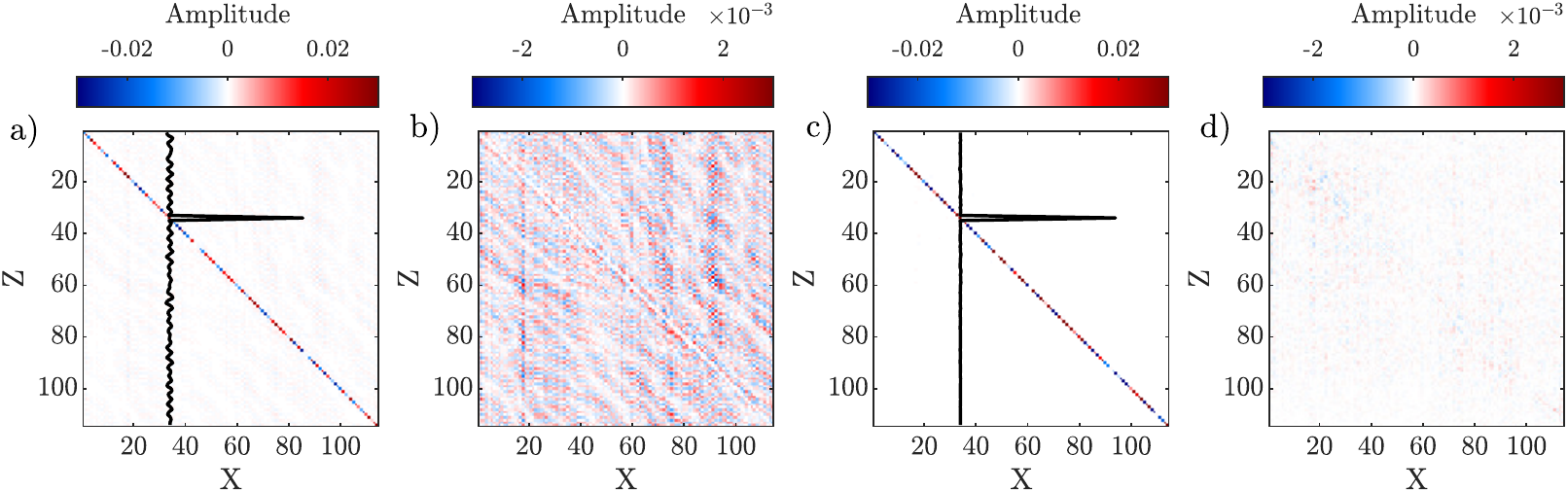}
\caption{The estimated real components of source signature matrix, i.e. $\bold{\Phi}^{\text{T}}\bold{A}_k\bold{U}^{k}$, (a-b) with the rough initial model 4(Fig. \ref{Fig_Mar}e) for 3 Hz, and (c-d) with the kinematically accurate velocity model 2(Fig. \ref{Fig_Mar}c) for 12 Hz. (a,c) shows $\bold{\Phi}^{\text{T}}\bold{A}_k\bold{U}^{k}$ where the column at $X=34$ of these matrices is plotted in black to have a visual insight about the magnitude of diagonal and off-diagonal elements. (b,d) show the non-diagonal components of (a,c) [$\bold{\Phi}^{\text{T}}\bold{A}_k\bold{U}^{k}-\text{Diag}(\text{diag}(\bold{\Phi}^{\text{T}}\bold{A}_k\bold{U}^{k}))$].}
\label{Fig_Mar_MAT}
\end{figure}

\begin{figure}
\centering
\includegraphics[scale=1]{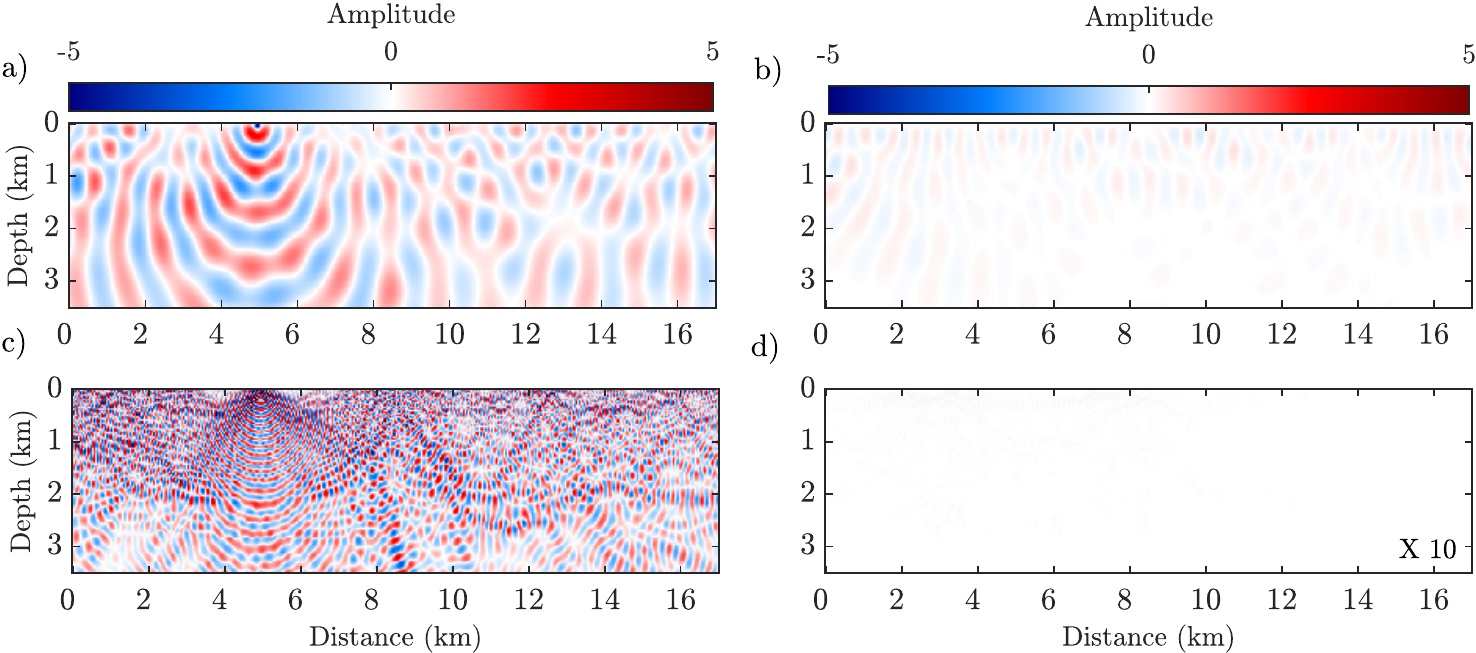}
\caption{(a) The reconstructed wavefield using Eq. \eqref{U} using the rough initial model 4 (Fig. \ref{Fig_Mar}e) for 3 Hz. (b) The difference between (a) and the reconstructed wavefield with known true source signature. (c-d) Same as (a-b), but using the kinematically accurate velocity model 2 (Fig. \ref{Fig_Mar}c) for 12~Hz.}
\label{Fig_Mar_Wavefield}
\end{figure}

\begin{figure}
\centering
\includegraphics[scale=1.0]{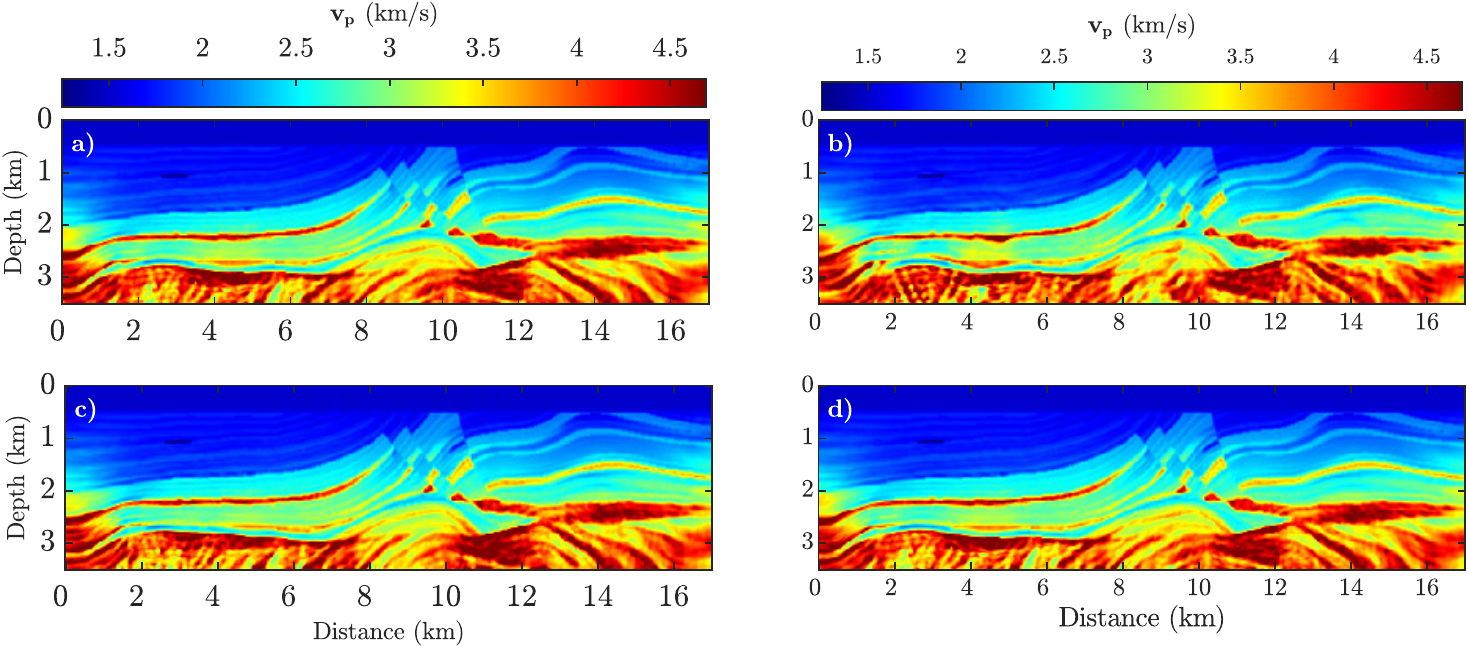}
\caption{Marmousi II inversion results when the inversion starts from the crude initial velocity model 4 (Fig. \ref{Fig_Mar}e). Estimated velocity model using IR-WRI (a) with known sources, (b-d) with unknown sources using (b) separate method, (c) joint method with Algorithm \ref{alg1}, and (d) joint method with Algorithm \ref{alg2}.}
\label{Fig_Mar_inv}
\end{figure}

\begin{figure}
\centering
\includegraphics[scale=1.0]{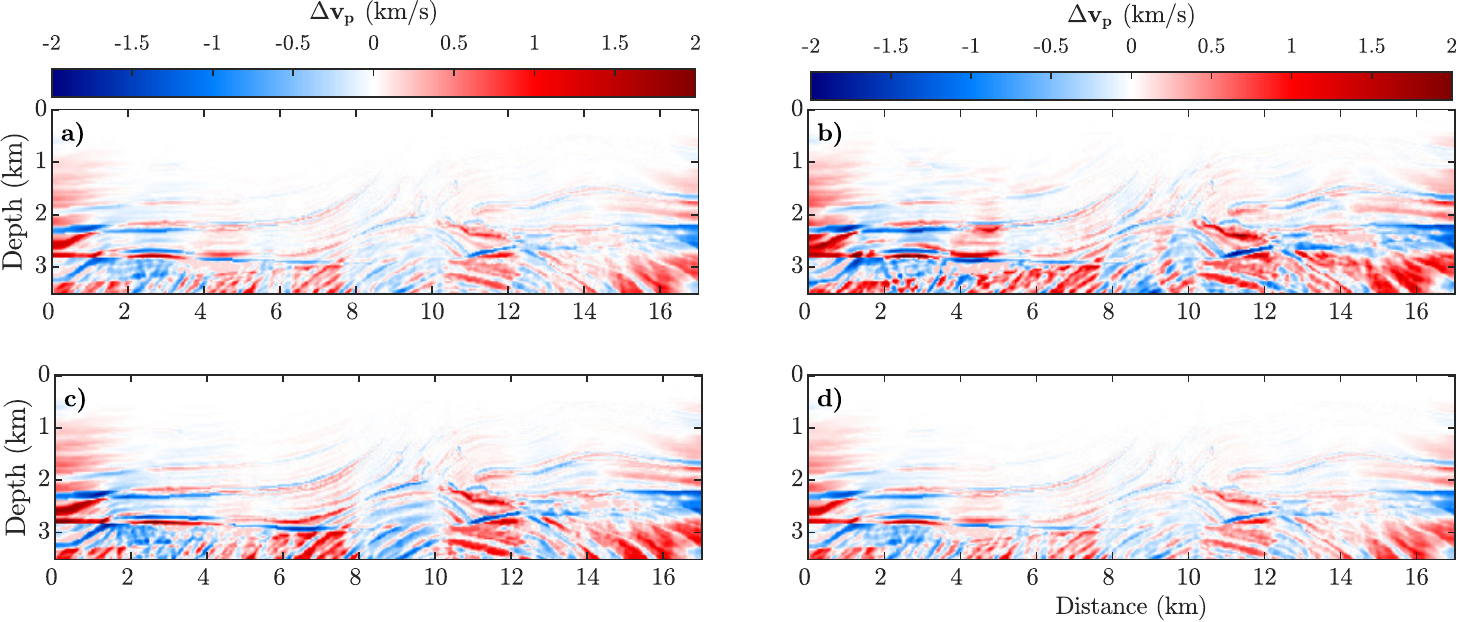}
\caption{Marmousi II inversion results. (a-d) The model error (difference between true model and estimated model) for estimated models Figs. \ref{Fig_Mar_inv}a-\ref{Fig_Mar_inv}d. }
\label{Fig_Mar_inv_res}
\end{figure}

\begin{figure}
\centering
\includegraphics[scale=0.7]{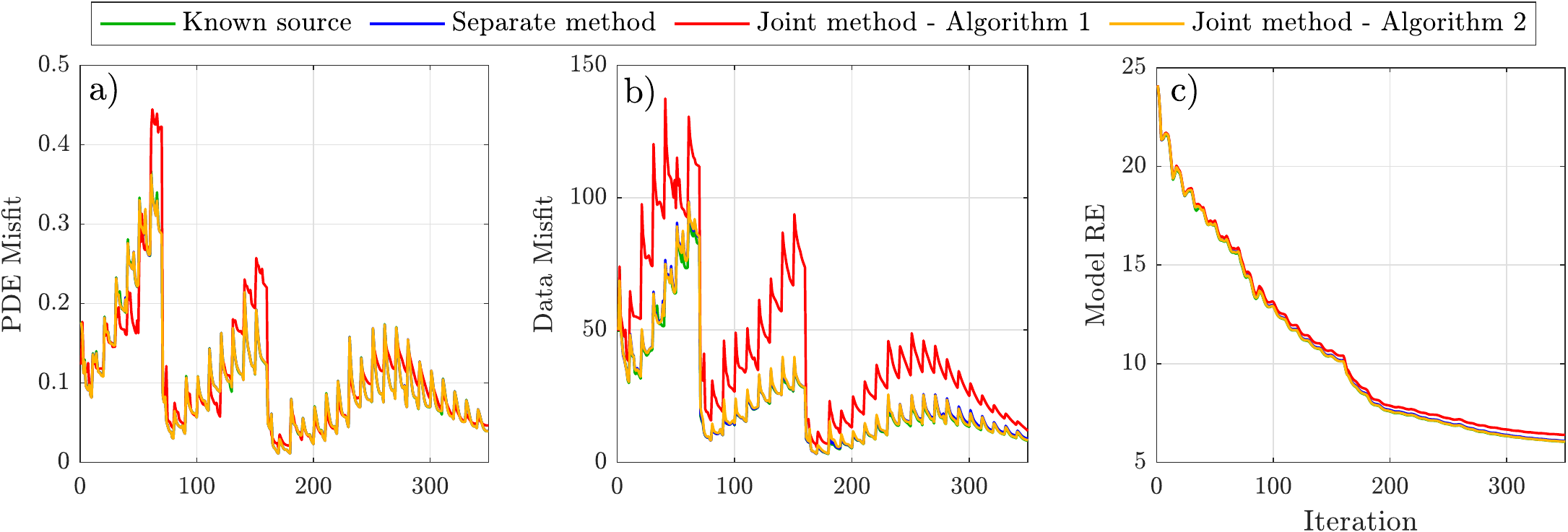}
\caption{Marmousi II test. A comparison between the convergence history of IR-WRI with known source (green), unknown sources with separate method (blue), joint method with Algorithm 1 (red) and joint method with Algorithm 2 (orange) as a function of iteration number. (a) PDE misfit, (b) Data misfit, and (c) Model RE.   }
\label{Fig_Mar_res}
\end{figure}

We continue by assessing Algorithms \ref{alg1} and \ref{alg2} as well as the separate method for FWI on the Marmousi II model when the inversion is started using the rough initial model 4 (Fig. \ref{Fig_Mar}e) and a 3~Hz frequency. 
The fixed-spread surface acquisition consists of 114 sources spaced 150~m apart with the source signatures depicted in Fig. \ref{Fig_wavelet} at the surface, and 340 hydrophone receivers spaced 50~m apart at 75~m depth.\\
We first show the effects of the source blending in the matrix $\bold{S}$, Eq. \eqref{S}, for the rough initial model 4 (Fig. \ref{Fig_Mar_MAT}a). The column of the matrix at $X=34$ is plotted in black to give more precise insights on the relative magnitude of the diagonal and off-diagonal elements. Also, the off-diagonal elements of this matrix are plotted separately in Fig. \ref{Fig_Mar_MAT}b. It is shown that the diagonal elements of this matrix are dominant, the maximum amplitude of the off-diagonal components being less than 1 percent of the maximum diagonal element. We remind that such off-diagonal elements partially absorb the errors in the estimated physical source signatures when the initial velocity model is rough, but we need to remove these effects during the inversion, which is the goal of Algorithms \ref{alg1} and \ref{alg2}.\\ 
The effects of the source blending are also seen in the reconstructed wavefields. The reconstructed monochromatic wavefield associated with the source located in Fig. \ref{Fig_Mar_MAT}a is shown in Fig. \ref{Fig_Mar_Wavefield}a. To assess the effects of the source blending, we show the difference between this wavefield and the extended wavefield reconstructed with the true source signature in Fig. \ref{Fig_Mar_Wavefield}b. It can be seen that the differences are not significant. 
This is the worst scenario because of the significant inaccuracy of the velocity model 4. However, the estimated source signature matrix becomes close to a diagonal matrix as soon as the accuracy of the velocity model improves. This statement is verified in Figures \ref{Fig_Mar_MAT}c-\ref{Fig_Mar_MAT}d and \ref{Fig_Mar_Wavefield}c-\ref{Fig_Mar_Wavefield}d, which are similar to Figures \ref{Fig_Mar_MAT}a-\ref{Fig_Mar_MAT}b and  \ref{Fig_Mar_Wavefield}a-\ref{Fig_Mar_Wavefield}b, except that the source signature and the wavefield are now estimated from the kinematically accurate velocity model 2 and a frequency of 12~Hz. It can be seen that the source matrix, Fig. \ref{Fig_Mar_MAT}c, tends to a diagonal matrix, and the differences between the estimated and the true wavefields tend to zero in Fig. \ref{Fig_Mar_Wavefield}d.\\
We continue by performing the frequency-domain IR-WRI in the 3~Hz - 12~Hz frequency band with a frequency interval of 0.5~Hz. Mono-frequency batches are successively inverted following a classical frequency continuation strategy. We perform three paths through the frequency batches to improve the inversion results, using one path's final model as the initial model of the next one. The starting and finishing frequencies of the paths are [3, 6], [3, 7], [3, 12]~Hz.
We compare the results of IR-WRI with the known sources \citep{Aghamiry_2019_IWR} (Fig. \ref{Fig_Mar_inv}a), IR-WRI with unknown sources using the separate method (Fig. \ref{Fig_Mar_inv}b), the joint method with Algorithm \ref{alg1} (Fig. \ref{Fig_Mar_inv}c), and finally, the joint method with Algorithm \ref{alg2} (Fig. \ref{Fig_Mar_inv}d), to assess the effectiveness of the proposed methods. Also, the model error (the difference between the estimated and true velocity model) for the different estimated models (Figs. \ref{Fig_Mar_inv}a-\ref{Fig_Mar_inv}d) are shown in Fig. \ref{Fig_Mar_inv_res}a-\ref{Fig_Mar_inv_res}d.
It can be seen that all of these methods work well but with a different computational cost. Also, IR-WRI with known sources (Fig. \ref{Fig_Mar_inv}a) and IR-WRI using Algorithm \ref{alg2} (Fig. \ref{Fig_Mar_inv}d) are really close together and have a better performance at the reservoir level (compare Figs. \ref{Fig_Mar_inv_res}a and \ref{Fig_Mar_inv_res}d). \\
Finally, the wave-equation residual, data residual, and RE for IR-WRI with known and unknown sources are shown in Fig. \ref{Fig_Mar_res}.  
In summary, we see that the separate method and both of the joint-approach algorithms work well and can reconstruct velocity models close to what we get from IR-WRI with known source signatures but with a different computational burden. The test in the next section shows that this conclusion is not valid for a more complicated velocity model.  
\subsection{2004 BP salt model}
\begin{figure}
\centering
\includegraphics[scale=0.5]{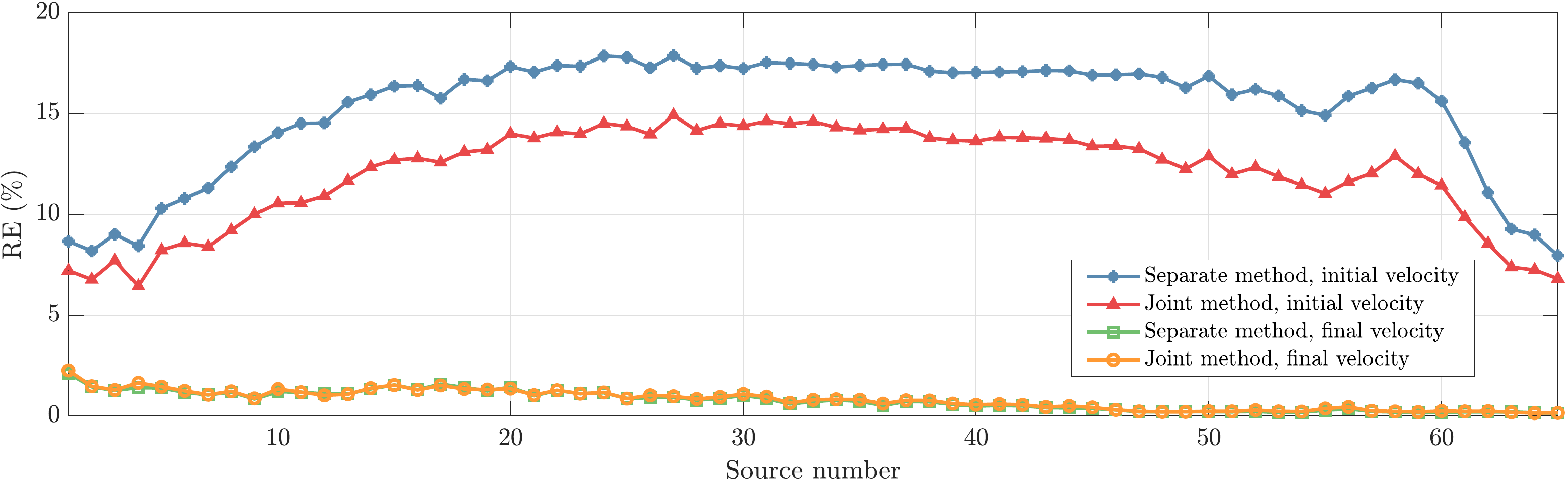}
\caption{BP Salt test. RE of the estimated source signatures for separate and the joint method when a rough initial model (Fig. \ref{Fig_BP_inv}b) and the accurate velocity model (Fig. \ref{Fig_BP_inv}f) are used.}
\label{Fig_BP_wavelet}
\end{figure}

\begin{figure}
\centering
\includegraphics[scale=1.2]{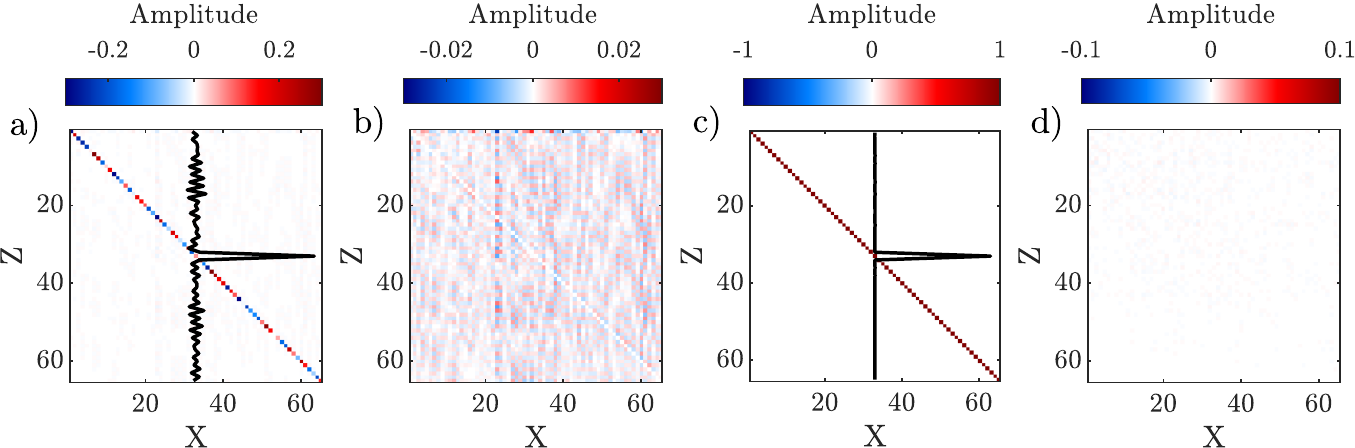}
\caption{BP Salt test. The estimated real components of source signature matrix, i.e. $\bold{\Phi}^{\text{T}}\bold{A}_k\bold{U}^{k}$, for the (a-b) first and (c-d) final iteration of IR-WRI. (a,c) shows $\bold{\Phi}^{\text{T}}\bold{A}_k\bold{U}^{k}$ where the column at $X=34$ of these matrices is plotted in black to have a visual insight about the magnitude of diagonal and off-diagonal elements. (b,d) show the non-diagonal components of (a,c) [$\bold{\Phi}^{\text{T}}\bold{A}_k\bold{U}^{k}-\text{Diag}(\text{diag}(\bold{\Phi}^{\text{T}}\bold{A}_k\bold{U}^{k}))$].}
\label{Fig_BP_MAT}
\end{figure}

\begin{figure}[htb!]
\centering
\includegraphics[scale=1.3]{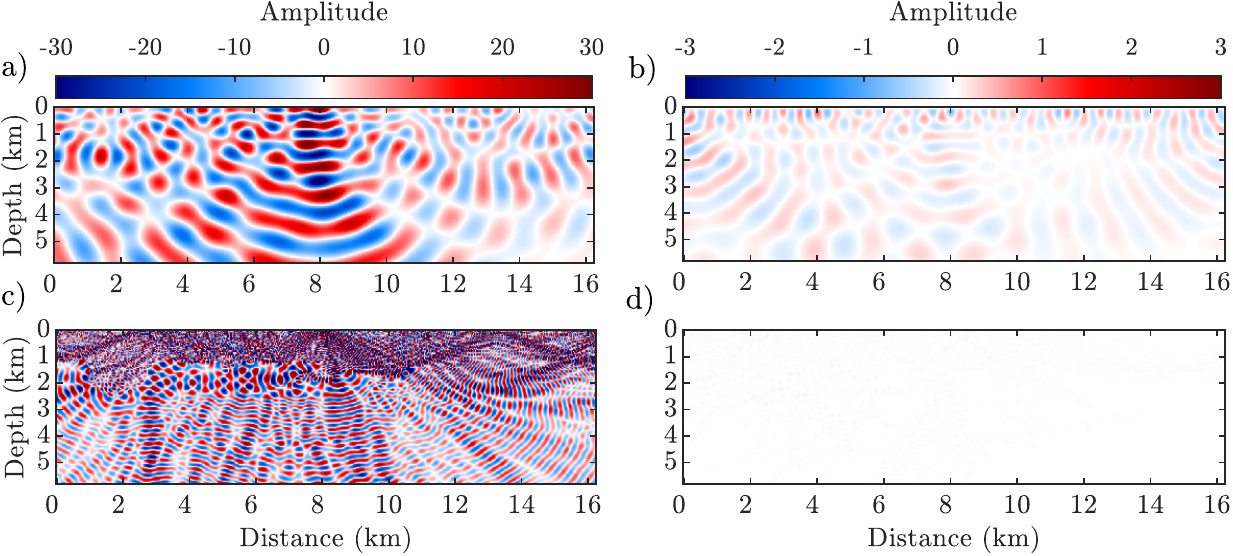}
\caption{BP Salt test. (a) The reconstructed wavefield using Eq. \eqref{U} at the first iteration of IR-WRI using joint method with the rough initial model (Fig. \ref{Fig_BP_inv}b) for 3 Hz. (b) The difference between (a) and the reconstructed wavefield with true source. (c-d) Same as (a-b), but for 12~Hz at the final iteration of IR-WRI with updated velocity model Fig. \ref{Fig_BP_inv}f.}
\label{Fig_BP_Wavefield}
\end{figure}

\begin{figure}[htb!]
\centering
\includegraphics[scale=1.3]{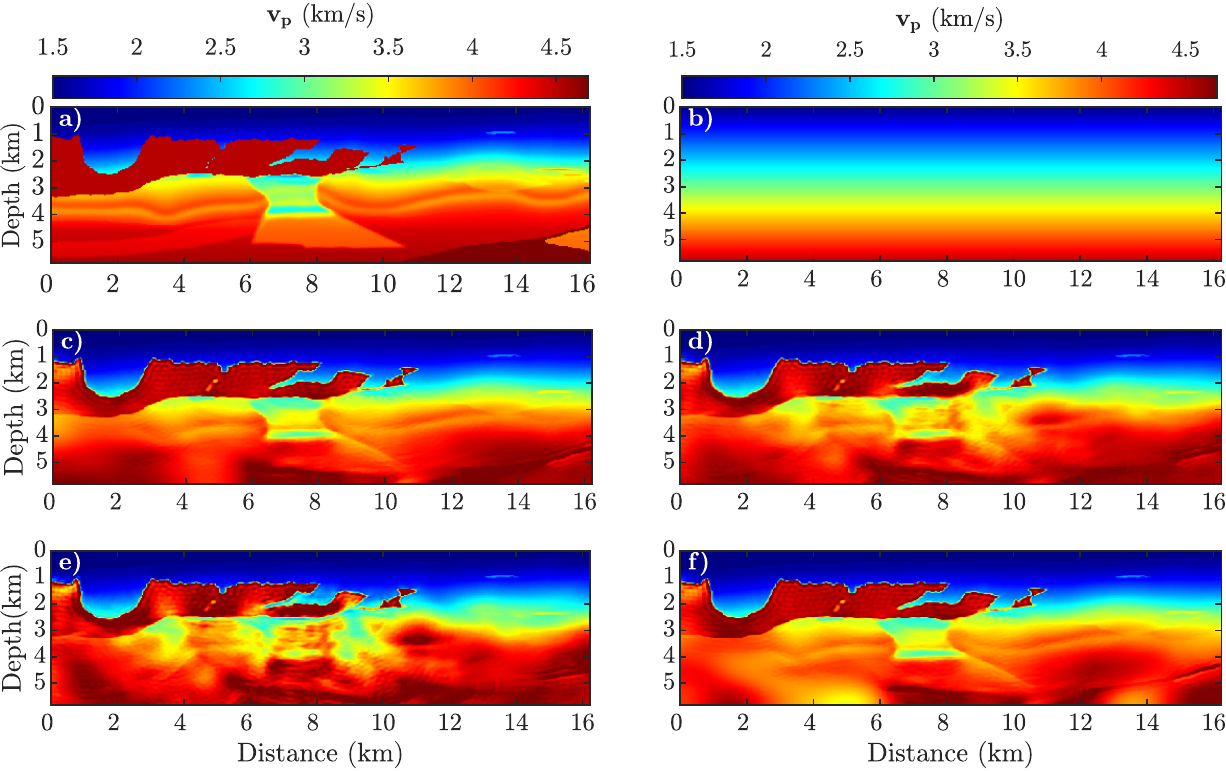}
\caption{2004 BP salt test: (a) True velocity model. (b) Initial velocity model. (c-f) BTT-regularized IR-WRI (c) with known sources, (d) with unknown sources using separate method, (e) with unknown sources using Algorithm \ref{alg1}, and (f) using Algorithm \ref{alg2}.}
\label{Fig_BP_inv}
\end{figure}

\begin{figure}[htb!]
\centering
\includegraphics[scale=1.3]{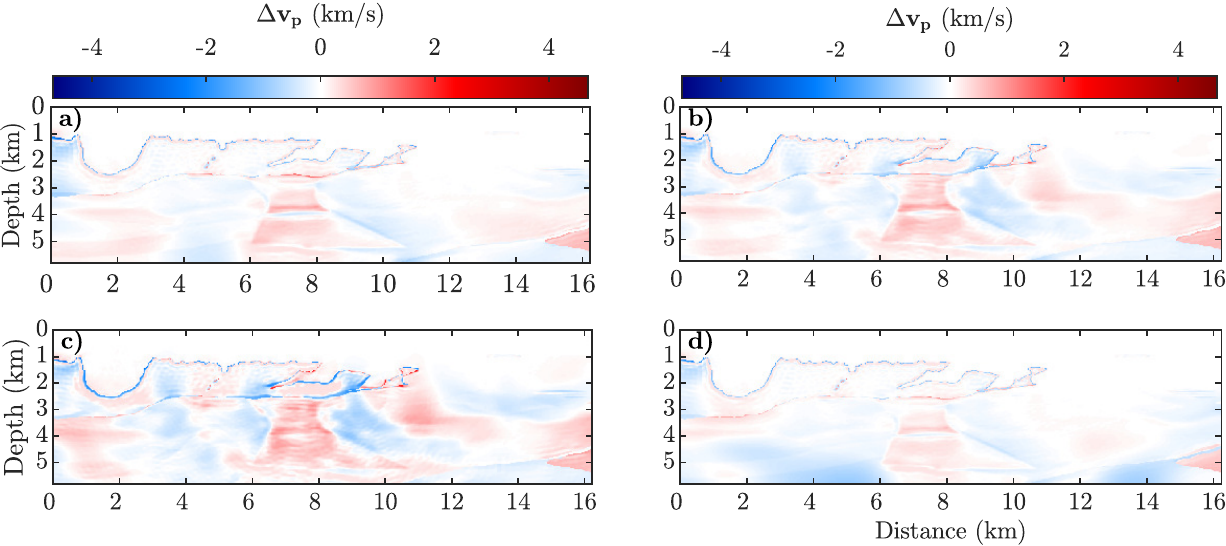}
\caption{2004 BP salt test: (a-d) The model error (difference between true model and estimated model) for estimated models Figs. \ref{Fig_BP_inv}c-\ref{Fig_BP_inv}f.}
\label{Fig_BP_inv_res}
\end{figure}

We continue by assessing the performance of the methods on a re-scaled version of left-target of challenging 2004 BP salt model \citep{Billette_2004_BPB} (Fig. \ref{Fig_BP_inv}a) when a 1D gradient initial model (Fig. \ref{Fig_BP_inv}b) and the 3~Hz frequency are used to start the inversion.
We use 65 point sources with 250~m spacing and a line of receivers with 50~m spacing at 75~m depth. Like the previous test, the source signatures are random Ricker wavelets with different central frequencies between [8-12]~Hz and the initial phases between [0-0.4]~s.
We apply the inversion in the 3~Hz-13~Hz frequency band with a frequency interval of 0.5~Hz. We perform three paths through the frequency batches, using the final model of one path as the initial model of the next one, and each batch contains two frequencies with one frequency overlap. The starting and finishing frequencies of the three paths are [3, 6], [4, 8.5], [6, 13]~Hz, respectively. 
Bound-constrained Tikhonov + Total variation (BTT)-regularization \citep{Aghamiry_2019_CRO} is applied on IR-WRI for all the cases to decrease the ill-posedness of the problem.
We first plot the RE for the estimated source signatures of the separate and joint method as a function of the source number in Fig. \ref{Fig_BP_wavelet}  for the first (3~Hz, using initial velocity model Fig. \ref{Fig_BP_inv}b) and the final iteration of the inversion (13~Hz, using updated velocity model \ref{Fig_BP_inv}f).
We see the same effects as those revealed by Marmousi II in Fig. \ref{Fig_err_wavelet} in the sense that the joint approach provides more accurate source signature estimation. 
%
%
%
Also, the estimated source signature matrices, Eq. \eqref{S}, at the first iteration and the final iteration of IR-WRI, are shown in Fig. \ref{Fig_BP_MAT} and the related reconstructed monochromatic wavefields are shown in Fig. \ref{Fig_BP_Wavefield}. Again, we see similar effects as those revealed by the Marmousi II test, except that the off-diagonal elements are larger relative to the diagonal counterparts, hence revealing the more complex structure of the BP salt model. 
We show in Fig. \ref{Fig_BP_inv} the final results of BTT-regularized IR-WRI with the known sources (Fig. \ref{Fig_BP_inv}c), the unknown sources using the separate method (Fig. \ref{Fig_BP_inv}d), the joint method with Algorithm 1 (Fig. \ref{Fig_BP_inv}e), and finally, the joint method with Algorithm 2 (Fig. \ref{Fig_BP_inv}f).
Also, the model errors for different estimated models are shown in Fig. \ref{Fig_BP_inv_res}.
In contrast to the Marmousi II test, the different methods don't converge to the same minimizer. Let's consider the BTT-regularized IR-WRI with known sources as the benchmark model (Fig. \ref{Fig_BP_inv}c). Only IR-WRI with unknown sources using Algorithm \ref{alg2} reaches approximately the same results. The failure of the separate method (Fig. \ref{Fig_BP_inv}d and \ref{Fig_BP_inv_res}b) probably results from the limited quality of the estimated wavelets at the early iterations when the initial velocity model is inaccurate (Fig. \ref{Fig_BP_wavelet}). Also, the failure of Algorithm \ref{alg1} (Fig. \ref{Fig_BP_inv}e and \ref{Fig_BP_inv_res}c) may result from the significant amplitudes of the off-diagonal elements of the estimated source signature matrix (Fig. \ref{Fig_BP_MAT}), and it seems that the iterative refinement implemented in Algorithm \ref{alg1} is not enough to correct all of these effects. For such complicated velocity models, we need to recompute the wavefields from the diagonalized source signature with Algorithm \ref{alg2} (lines 4-5).
\section{Conclusions}
We extended the recently proposed iteratively-refined wavefield reconstruction inversion (IR-WRI) to estimate the unknown source signatures.
The source signatures and wavefields are jointly estimated with a variable projection method during the extended wavefield reconstruction subproblem.
We first show that the source signature estimation generates computational overhead when each source is processed separately because the extended wave equation operator becomes source dependent. This computational overhead becomes prohibitive when the augmented wave equation system is solved with the direct method as a LU factorization needs to be performed for each source.   To bypass this issue and make the operator source independent, we proposed a method that blends the sources during each wavefield reconstruction. Accordingly, for each source of the experiment, the proposed method searches for a virtual blended source that best fits each single-source dataset during the wavefield reconstruction. Regardless of the computational efficiency issue, we also show that when the background velocity model is inaccurate, this source blending provides a more accurate source signature estimation at the physical source location than the case where source blending is not used. This probably results from the additional degrees of freedom provided by the extra virtual sources.
Once the spatially distributed source signatures have been estimated, we restrict them at the position of the physical source to mimic the true source signatures, and we re-estimate the extended wavefields with these localized source signatures. Although we assume in this study that the source locations match the grid points of the computational domain, the method can be readily used for arbitrary source positions.
\section{ACKNOWLEDGEMENTS}  
This study was partially funded by the WIND consortium (\href{https://www.geoazur.fr/WIND}{https://www.geoazur.fr/WIND}), sponsored by Chevron, Shell, and Total. This study was granted access to the HPC resources of SIGAMM infrastructure (\href{http://crimson.oca.eu}{http://crimson.oca.eu}), hosted by Observatoire de la C\^ote d'Azur and which is supported by the Provence-Alpes C\^ote d'Azur region, and the HPC resources of CINES/IDRIS/TGCC under the allocation A0050410596 made by GENCI. 

\bibliographystyle{SEG}
\newcommand{\SortNoop}[1]{}

\end{document}